\documentclass[11pt]{amsart}
\usepackage{amssymb,latexsym,amsmath,graphicx,amsthm}

\makeatletter
\@addtoreset{figure}{section}
\def\thefigure{\thesection.\@arabic\c@figure}
\def\fps@figure{h, t}
\@addtoreset{table}{bsection}
\def\thetable{\thesection.\@arabic\c@table}
\def\fps@table{h, t}
\@addtoreset{equation}{section}

\makeatother

\def\grad{\boldsymbol{\nabla}}

\def\NR{\nonumber \\} 
\newtheorem{thm}{Theorem}[section]
\newtheorem{prop}[thm]{Proposition}
\newtheorem{lem}[thm]{Lemma}

\newtheorem{dfn}[thm]{Definition}
\newtheorem{cnj}{Conjecture}

\newcommand{\Maps}{\operatorname{Maps}}

\newcommand{\CC}{{\mathbb C}}
\newcommand{\RR}{{\mathbb R}}
\renewcommand{\SS}{{\mathbb S}}
\newcommand{\ZZ}{{\mathbb Z}}

\newcommand{\fg}{{\mathfrak g}}
\newcommand{\fh}{{\mathfrak h}}
\newcommand{\fp}{{\mathfrak p}}

\renewcommand{\P}{{\mathcal P}}
\newcommand{\K}{{\mathcal K}}

\begin{document}

\title{Combinatorial formulas for products of Thom classes}

\author[V. Guillemin]{Victor Guillemin\footnotemark {*}}
\thanks{* Supported by NSF grant DMS 890771.} 
\address{Department of Mathematics, MIT, Cambridge, MA 02139}
\email{vwg@@math.mit.edu}
\author[C. Zara]{Catalin Zara \footnotemark{**}}
\thanks{** This research was partially 
conducted by the second author for the Clay Mathematics Institute.}
\address{Department of Mathematics, MIT, Cambridge, MA 02139}
\email{czara@math.mit.edu}

\begin{abstract}
Let $G$ be a torus of dimension $n > 1$ and $M$ a compact 
Hamiltonian $G$-manifold with $M^G$ finite. A circle, $S^1$, in $G$ is 
\emph{generic} if $M^G = M^{S^1}$. For such a circle the moment map 
associated with its action on  $M$ is a perfect Morse function. 
Let $\{ W_p^+ ; p \in M^G\}$ be the Morse-Whitney stratification of 
$M$ associated with this function, and let $\tau_p^+$ be the equivariant 
Thom class dual to $W_p^+$. These classes form a basis of $H_G^*(M)$ as a 
module over $\SS(\fg^*)$ and, in particular,
$$\tau_p^+ \tau_q^+ = \sum c_{pq}^r \tau_r^+$$
with $c_{pq}^r \in \SS(\fg^*)$. For manifolds of GKM type we obtain a 
combinatorial description of these $\tau_p^+$'s and, from this description, 
a combinatorial formula for $c_{pq}^r$.
\end{abstract}

\maketitle

\section{Products of Thom classes}
\label{sec:1}

Let $M^{2d}$ be a compact Hamiltonian $S^1$-manifold with moment
map, $\phi :M \to \RR$.  If $M^{S^1}$ is finite, $\phi $ is a
Morse function, and its critical points, $p \in M^{S^1}$, are all 
of even index.  This has important consequences for the topology
of $M$:  If we put an $S^1$-invariant Riemannian metric, $B$, on $M$ and let
$$  v=\grad_B \phi $$
be the gradient vector field associated with $B$ and $\phi$,
then, for every critical point, $p \in M^{S^1}$, the unstable
manifold at $p$:
$$ W^+_p= \{ q \in M ; \,\, \lim_{t\to -\infty}
    (\exp tv) (q) =p \} $$
supports a cohomology class, and these ``Thom'' classes
\begin{equation}
  \label{eq:1.1}
  \tau^+_p, \quad p\in M^{S^1}
\end{equation}
are an additive basis of the cohomology ring, $H^* (M,\RR )$.
The same is true of the stable manifolds
$$  W^-_p =\{ q \in M ; \,\, \lim_{t \to \infty}
  (\exp tv ) (q) =p \} $$
and their dual Thom classes
$$\tau^-_p, \quad p \in M^{S^1}.$$
The main topic of this paper will be the symplectic version of
what is sometimes called the multiplicative Morse problem:  
Given $p$ and $q$ in $M^{S^1}$, 
$\tau^+_p \tau^+_q$ can be expanded as a sum
$$  \tau^+_p \tau^+_q = \sum c^r_{pq} \tau^+_r \, .$$
What are the $c^r_{pq}$'s?  Closely related to this is the
question of determining the cohomology pairings:
\begin{equation}
    \label{eq:1.2}
c_{pqr} = \int \tau^+_p \tau^+_q \tau^-_r \; .
\end{equation}
Neither of these questions is easy to answer even when the
structure of the cohomology ring itself is well understood.  For
instance if $M$ is the coadjoint orbit of a compact Lie group,
the computation of the $c^r_{pq}$'s is an important open problem
in the theory of the Schubert calculus and is the focus of a 
lot of recent activity.  
(See, for instance, [BGG], [BH], [Bi], [Ko] and  [Kn].)

In this paper we will consider the \emph{equivariant} version of
this problem.  We will assume  the action of $S^1$ on $M$ can 
be extended to a Hamiltonian action of a torus, $G$, of dimension 
$n>1$, and replace the $\tau^+_p$'s in \eqref{eq:1.2} by their
equivariant counterparts.  These equivariant Thom classes generate 
$H_G (M,\RR)$ as a module
over the ring, $H_G (pt)= \SS(\fg^*)$, so one gets as above an
identity
$$  \tau^+_p \tau^+_q = \sum c^r_{pq} \tau^+_r \; , $$
but now the $c^r_{pq}$'s are elements of the
polynomial ring, $\SS (\fg^*)$.  When 
$$
\text{degree }\tau^+_r= \text{degree }\tau^+_p + 
\text{degree }\tau^+_q \; ,  $$ 
they are polynomials of degree
zero (\emph{i.e.} real numbers) and, in fact, coincide with the $c^r_{pq}$'s
in \eqref{eq:1.2}.  Thus, they are in principle a much larger
list of unknown quantities.  We will show, however, that they are in some
sense, easier to compute due to the fact that, in equivariant
cohomology, one has a much richer store of intersection
invariants to play around with.  More explicitly if $X$ and $Y$ 
are submanifolds of $M$ and $\tau_X$ and $\tau_Y$ their dual Thom 
classes, the intersection number
\begin{equation}
  \label{eq:1.3}
  \# (X \cap Y) = \int \tau_X \tau_Y
\end{equation}
is zero except when $X$ and $Y$ are of complementary dimension.
On the other hand if $X$ and $Y$ are $G$-invariant and $\tau_X$
and $\tau_Y$ their equivariant Thom classes, the expression
\eqref{eq:1.3} (which is now an element of $\SS(\fg^*)$) can be
non-zero no matter what the relative dimensions of $X$ and $Y$
are.  Moreover, by the localization theorem of
Atiyah-Bott-Berline-Vergne, \eqref{eq:1.3}, is a sum of local
intersection invariants
$$  \# (X \cap Y)_p \in Q (\fg^*) \, , $$
where $Q(\fg^*)$ is the quotient field of $\SS(\fg^*)$ and $p$ a fixed 
point, and each of these is itself an intersection invariant.

Of particular interest for us will be certain intersection
invariants of this type associated with the moment map, $\phi$.
Suppose that $p$ and $q$ are critical points of $\phi$ and that
there are no critical values of $\phi$ in the interval 
$(\phi (p), \phi (q))$.  Let $\phi (p)<c< \phi (q)$ and let
$$ M_c = \phi^{-1} (c) / S^1 $$
be the symplectic reduction of $M$ at $c$.  By the
Marsden-Weinstein theorem, $M_c$ is a symplectic orbifold, and
the action of $G$ on $M$ induces an action of the group
$$G_1 = G/S^1$$
on $M_c$.  The reduced spaces
$(W^+_p)_c$ and $(W^-_q)_c$ 
are $G$-invariant symplectic sub-orbifolds of $M_c$ and so their
equivariant intersection ``number'' 
\begin{equation}
  \label{eq:1.4}
  \# ((W^+_p)_c \cap (W^-_q)_c)
\end{equation}
is well-defined as an element of the subring, $\SS(\fg^*_1)$, of $
\SS(\fg^*)$.  Moreover, if  $M_c^{G_1}$  is finite,
then for every $v \in M^{G_1}_c$,  the local
intersection number
\begin{equation}
  \label{eq:1.5}
  \# ((W^+_p)_c \cap (W^-_q)_c)_v
\end{equation}
is well-defined as an element of $Q(\fg^*_1)$.

We will now describe the role of these intersection numbers in
the computation of the $c^r_{pq}$'s.  We will say that $M$ is 
\emph{a GKM manifold} if, for all non-critical values, $c$, 
 $M_c^{G_1}$ is finite.  Thus, being GKM is a
necessary and sufficient condition for the invariants
\eqref{eq:1.5} to be well-defined.  We recall some other
characterizations of these manifolds.

\begin{thm}
  \label{th:1.1}
  $M$ is GKM if and only if, for every $p \in M^G$, the weights,
  $\alpha_{i,p}$, of the linear isotropy representation of $G$ on 
  $T_pM$ are pair-wise linearly independent; i.e,~$\alpha_{i,p} $ 
  is not a multiple of $\alpha_{j,p}$ if $i \neq j$.
\end{thm}

\begin{thm}
  \label{th:1.2}
  $M$ is GKM if and only if, for every codimension one subtorus, $ 
 H$, of $G$ the connected components of $M^H$ are either fixed
 points of $G$ or imbedded copies of $S^2$.  Moreover, if a
 connected component, $X$, is a copy of $S^2$, the action of $G$
 on $X$ is symplectomorphic to the standard action of $G/H=S^1$ on $S^2$.
\end{thm}

\begin{thm}
  \label{th:1.3}
  $M$ is GKM if and only if the one skeleton of $M$
  \begin{equation}
    \label{eq:1.6}
    \{ p \in M, \quad \dim G_p \geq n-1 \}
  \end{equation}
 is a finite union of embedded $S^2$'s.
\end{thm}

\begin{proof}
  Theorem~\ref{th:1.3} is an obvious consequence of
  Theorem~\ref{th:1.2}; and it is easy to see that if the
hypotheses of Theorem~\ref{th:1.2} hold, $M$ is GKM.  For the
other implications see \cite{GZ2}.
\end{proof}

The intersection properties of the embedded $S^2$'s in the set
\eqref{eq:1.6} can be described by an intersection graph, and 
a beautiful observation of Goresky-Kottwitz-MacPherson is that
one can read off the structure of the equivariant cohomology ring 
of $M$ from the ``action'' of $G$ on this graph.  More explicitly 
let $\Gamma$ be the graph whose vertices are the fixed points of
$G$ and whose edges, $e$, are copies, $X_e$, of the $S^2$'s in
Theorem~\ref{th:1.3}.  The graph structure on this collection of vertices
and edges is given by defining the pair of vertices incident to 
an edge, $e$, to be the set
$$  \partial e = \{ p,q \} = X^G_e \, . $$
In particular, if $e$ is an \emph{oriented} edge, its initial
vertex, $i(e)$, is defined to be the ``south pole'', $p$, of the
two-sphere, $X_e$, and its terminal vertex, $t(e)$, to be the
``north pole'', $q$ of $X_e$.  
  The action of $G$ on the set \eqref{eq:1.6} can be described
  graph theoretically by two pieces of data:  a function $\rho$
  which assigns to each oriented edge, $e$, of $\Gamma$ a
  one-dimensional representation, $\rho_e$, of $G$ and a function, 
$\kappa$, which assigns to each vertex, $p$, a $d$-dimensional 
 representation, $\kappa_p$.  These functions are defined by
 letting $\rho_e$ be the representation of $G$ on $T_p X_e$, $p =i(e)$,
 and letting $\kappa_p$ be the
 representation of $G$ on $T_pM$.  It is easily checked that 
$\rho$ and $\kappa$ satisfy the  axioms:  

\begin{align}
  \label{eq:1.7}
  \kappa_p &= \bigoplus_{i(e) =p} \rho_e \\
\label{eq:1.8}
\rho_{\bar{e}} &= \rho^*_e\\
\intertext{and}
\label{eq:1.9}
\kappa_p |_{G_e} &= \kappa_q |_{G_e}
\end{align}
where $\bar{e}$ is the edge obtained from $e$ by reversing
its orientation, $G_e$ is the kernel of $\rho_e$ and 
$p$ and $q$ are the vertices $i(e)$ and $t(e)$.  In 
particular, by \eqref{eq:1.7}, $\kappa_p$ is determined by the $
\rho_e$'s with $i(e)=p$.  Since $\rho_e$ is a one-dimensional
representation it is determined by its weight, $\alpha_e$; so the 
``action'' of $G$ on $\Gamma$ associated with $\rho$ and $\kappa$ 
consists essentially of a labeling of the oriented edges, $e$, of 
$\Gamma$ by weights, $\alpha_e$.  The axioms
\eqref{eq:1.7}--\eqref{eq:1.9} impose, of course, some
condition on this labeling.  For instance \eqref{eq:1.8} is
equivalent to
$$  \alpha_e =-\alpha_{\bar{e}} \, .$$

Now let $H_G(M)$ be the equivariant cohomology ring of $M$ and
$H_G (M^G)$ the equivariant cohomology ring of $M^G$.  Since
$M^G$ is a finite disjoint union of fixed points and these fixed
points are also the vertices, $V_{\Gamma}$, of $\Gamma$ it
follows that
$$  H_G (M^G) = \Maps (V_{\Gamma} , \SS(\fg^*)) \, . $$
Moreover, if $i : M^G \to M$ is the inclusion, the map 
$i^* : H_G(M) \to H_G(M^G)$ is injective, by a well-known result of Kirwan. 
The theorem of Goresky-Kottwitz-MacPherson which we alluded to
above asserts:

\begin{thm}
  \label{th:1.4}
A map $h: V_{\Gamma} \to S (\fg^*)$ is in 
the image of $i^*$ if and only if, for every edge, $e$, of
$\Gamma$
\begin{equation}
  \label{eq:1.10}
  h_p |_{\fg_e} = h_q |_{\fg_e}
\end{equation}
$p$ and $q$ being the vertices of $e$ and $\fg_e$ the annihilator 
of $\alpha_e$ in $\fg$.
\end{thm}

This leads us to define the \emph{equivariant cohomology ring}, 
$H (\Gamma , \alpha)$ of $\Gamma$ to be the set of all maps, 
$h:V_{\Gamma} \to S (\fg^*)$ satisfying \eqref{eq:1.10}.  Each of
the Thom classes \eqref{eq:1.1} gets mapped by $i^*$ onto an
element of $H (\Gamma , \alpha)$ and we will continue to use the
notation, $\tau^+_p$, for this ``combinatorial'' Thom class.  The 
main result of this article is a formula for this Thom class as a 
kind of path integral over certain paths in $\Gamma$.  Before
stating this result we'll describe a few basic
properties of these combinatorial Thom classes.  Lets continue to 
denote by $\phi $ the
restriction of the moment map, $\phi$, to $M^G$.  Identifying
$M^G$ with $V_{\Gamma}$, one can think of $\phi$ as a real-valued
function on $V_{\Gamma}$.  By Theorem~\ref{th:1.2}, $\phi$ takes on
distinct values on the vertices $i(e)$ and $t(e)$ of an oriented
edge, $e$.  We will say that this edge is \emph{ascending}
if $\phi (i(e))< \phi (t(e))$ and \emph{descending} if the
reverse inequality is true.  More generally if $\gamma$ is a path 
in $\Gamma$ we will say that $\gamma$ is \emph{ascending}
if each of its edges is ascending.  For every vertex, $p
\in V_{\Gamma}$ define the \emph{index}, $\sigma_p$, of $p$ to be 
the number of descending edges, $e$, with $i(e)=p$.

\begin{thm}
  \label{th:1.5}
The Thom class, $\tau^+_p$, has the following properties:

\begin{enumerate}
\item 
  Its support is the set of all vertices of $\Gamma$ which can be 
  joined to $p$ by an ascending path.

\item 
  The value of $\tau^+_p$ at $p$ is
$$  \nu_p^+ =   \sideset{}{'}\prod_{i(e)=p} \alpha_e \quad , $$
the product, $\prod'$, being over all descending edges with
$i(e)=p$. 
\end{enumerate}
\end{thm}

In certain instances these properties uniquely characterize $\tau^+_p$.

\begin{thm}
  \label{th:1.6}
Suppose that the indexing function, $\sigma : V_{\Gamma} \to
\ZZ$, $p \to \sigma_p$, is strictly increasing along
ascending paths.  Then $\tau^+_p$ is the \emph{unique}
element of $H(\Gamma, \alpha)$ with properties 1 and 2 above.
\end{thm}

We will now describe our ``path-integral'' formula for
$\tau^+_p$.  As mentioned above this formula will involve the
Hamiltonian action of the subgroup, $S^1$, of $G$ on $M$ and the
intersection invariants \eqref{eq:1.4} and \eqref{eq:1.5}.  Let
$\xi \in \fg$ be the infinitesimal generator of this subgroup and 
let $e$ be an ascending edge of $\Gamma$ with $p=i(e)$. 
For any point, $c$, on the interval between $\phi(p)$ 
and $\phi(q)$,  the $S^1$-reduced space, $(X_e)_c$ consists of
a single point, $v \in M_c$.  Let $\iota_e$ be the local
intersection number \eqref{eq:1.5}.

\begin{thm}
  \label{th:1.7}
For every $q \in V_{\Gamma}$
\begin{equation}
  \label{eq:1.11}
  \tau^+_p (q) = \sum E (\gamma)
\end{equation}
the sum being over all ascending paths, $\gamma$, joining
$p$ to $q$, and the summands being defined by
\begin{equation}
  \label{eq:1.12}
  E(\gamma) = (-1)^m \nu_q^+ \frac{\iota_{e_1}}{\hat{\alpha}_m}
  \prod^m_{k=2} \frac{\iota_{e_k}}{\hat{\alpha}_{k-1}-\hat{\alpha}_k} \; ,  
\end{equation}
where $e_1, \ldots ,e_m$ are the edges of $\gamma$ and 
$$  \hat{\alpha}_k = \frac{\alpha_{e_k}}{\alpha_{e_k} (\xi)} \, .$$
\end{thm}

\smallskip

\noindent {\bf Remarks:}

\begin{enumerate}
\item 
  The local intersection number, $\iota_e$, is equal to the
  global intersection number \eqref{eq:1.4} provided that
  there are no ascending paths in $\Gamma$ of length
  greater than one joining $p= i(e)$ to $q=t(e)$.  In particular, 
  if $\gamma$ is a longest path joining $p$ to $q$ all the
  intersection numbers in \eqref{eq:1.12} are global
  intersection numbers and in particular are elements of 
        $\SS(\fg^*)$.

\item  In Section~\ref{sec:4} we will give a purely
  combinatorial definition of $\iota_e$.  
\end{enumerate}

As a corollary of
  Theorem~\ref{th:1.5} one gets for \eqref{eq:1.2} the formula
  \begin{equation}
    \label{eq:1.13}
    c_{pqr} = \sum_t \delta_t E(\gamma_1) E(\gamma_2) E(\gamma_3) 
  \end{equation}
summed over all configurations of paths, $\gamma_1$, $\gamma_2$
and $\gamma_3$, where $\gamma_1$ is an ascending path from $p$ to
$t$, $\gamma_2$ an ascending path from $q$ to $t$, $\gamma_3$ a
descending path from $r$ to $t$ and
$$  \delta_t = \Bigl( \prod_{i(e)=t} \alpha_e \Bigr)^{-1} \, .$$

In particular:

\begin{thm}
 \label{th:1.8}
If the hypotheses of Theorem~\ref{th:1.6} are satisfied, 
$c^r_{pq} = c_{pqr}$ and hence $c^r_{pq}$ is equal to the
sum \eqref{eq:1.13}.

\end{thm}

A few words about the organization of this paper.  In
Section \ref{sec:2} we will give a brief account of the theory of 
$G$-actions on graphs 
(based, for the most part, on material in \cite{GZ3}).
In Section~\ref{sec:3} we will derive a preliminary version of
Theorem~\ref{th:1.7} and then, in Section ~\ref{sec:4}, deduce from it
the version above, after first describing how to define the
invariants \eqref{eq:1.5} combinatorially.  In Section~\ref{sec:5} we
will attempt to demystify what is perhaps the most puzzling feature of
the formula \eqref{eq:1.11}, the fact that all the summands are
rational functions (elements of the quotient field, $Q (\fg^*)$), 
whereas the sum itself is a polynomial.  This indicates that a
lot of mysterious cancellations are occurring in this summation;
and we will show \emph{how} these cancellations occur in a few simple 
but enlightening examples.

We would like to thank Tara Holm and Sue Tolman for helping us with the
computations involved in these examples and Allen Knutson for pointing out
to us antecedents in the combinatorial literature for formulas of type 
\eqref{eq:1.11} and \eqref{eq:1.13}.

\section{$G$-actions on graphs}\label{sec:2}

Let $\Gamma$ be a finite $d$-valent graph.  Given an oriented
edge, $e$, of $\Gamma$ we will denote by $i(e)$ the initial
vertex of $e$ and by $t(e)$ the terminal vertex; and we will
denote by $\bar{e}$ the edge obtained from $e$ by reversing its
orientation.  Thus $i(\bar{e})=t(e)$ and $t(\bar{e})=i(e)$.

\begin{dfn}
  \label{dfn:2.1}
Let $\rho$ be a function which assigns to each oriented edge,
$e$, of $\Gamma$ a one dimensional representation, 
$\rho_e: G \to S^1$; and let $\kappa$ be a function which assigns to each
vertex, $p$, a $d$-dimensional representation of $G, \kappa_p$.
We will say that $\rho$ and $\kappa$ define \emph{an action of 
$G$  on $\Gamma$} if the axioms \eqref{eq:1.7}--\eqref{eq:1.9}
are satisfied.  
\end{dfn}

Let $\alpha_e$ be the weight of the
representation $\rho_e$.  By \eqref{eq:1.7}, the weights,
$\alpha_e, i(e) =p$, determine the representation $\kappa_p$ up
to isomorphism; so the action of $G$ on $\Gamma$ is basically
just a labeling of the edges of $\Gamma$ by weights.  We will
denote the function, $e \to \alpha_e$, by $\alpha$ and call it
the \emph{axial function} of the action of $G$ on $\Gamma$.  The
axioms \eqref{eq:1.7}--\eqref{eq:1.9} can be reformulated as
statements about $\alpha$:

\begin{prop}
  \label{prop:2.2}
Axiom \eqref{eq:1.8} is satisfied iff
$$  \alpha_e = -\alpha_{\bar{e}}$$
and axiom \eqref{eq:1.9} is satisfied iff one can order the
weights
\begin{align}
\alpha_{e_k} \; , & \quad  i(e_k)=p, \, e_k \neq e \NR
\intertext{and the weights}
 \alpha_{e'_k}, \; , & \quad  i(e'_k)=q, \, e'_k \neq \bar{e} \nonumber
\end{align}
so that 
\begin{equation}
  \label{eq:2.1}
  \alpha_{e'_k} = \alpha_{e_k} + c_k \alpha_e \, .
\end{equation}

\end{prop}

(We will leave the proof of these assertions as an easy exercise.)

\begin{dfn}
  \label{dfn:2.3}
The action of $G$ on $\Gamma$ is a \emph{GKM action} if, for all
vertices, $p$, the weights, $\alpha_e, i(e)=p$, are pair-wise
linear independent.
\end{dfn}

From now on we will assume, unless we state otherwise, that the
action of $G$ on $\Gamma$ has this property.

For every vertex, $p$, of $\Gamma$ let $E_p$ be the set of
oriented edges, $e$, with $i(e)=p$.

\begin{dfn}
  \label{dfn:2.4}
A \emph{connection} on $\Gamma$ is a function which assigns to
every oriented edge, $e$, a bijective map
$$  \theta_e : E_{i(e)} \to E_{t(e)}$$
satisfying $\theta_{\bar{e}} = \theta^{-1}_e$.  This connection
is \emph{compatible with the action of } $G$ if, for every
oriented edge, $e$, with $i(e)=p$ and every edge, $e_{k} \in 
E_p, \, e_{k} \neq e$
\begin{equation}
  \label{eq:2.2}
  \alpha_{e'_k} =  \alpha_{e_k} +
  c_k \alpha_e, \quad \text{ where  } e'_k =\theta_e (e_k) \, .
\end{equation}

\end{dfn}

Thus the existence of a $G$-compatible connection is a slight
sharpening of the identity \eqref{eq:2.1}.  It is easy to see
that $G$-compatible connections exist, and we will assume
henceforth that $\Gamma$ is equipped with such a connection.

Let $V_{\Gamma}$ be the set of vertices of $\Gamma$ and
$E_{\Gamma}$ the set of oriented edges.  Motivated by the theorem 
of Goresky--Kottwitz--MacPherson we define the equivariant
cohomology ring, $H(\Gamma , \alpha)$, of $\Gamma$ to be the set
of maps, $h:V_{\Gamma} \to \SS (\fg^*)$, satisfying the
compatibility condition \eqref{eq:1.10} for all $e \in
E_{\Gamma}$.  
This ring has a natural 
grading\footnote{This definition is, unfortunately, 
inconsistent with the topological definition which assigns to  
$H^k$ the degree $2k$.  
(It is, however, more natural in this algebraic context.)} 
$$H^k (\Gamma , \alpha) = H(\Gamma, \alpha)
\cap \Maps (V_{\Gamma}, \SS^k (\fg^*))$$
and contains $\SS(\fg^*)$ as a subring:  the ring of constant maps
of $V_{\Gamma}$ into $\SS(\fg^*)$.  The proof of Theorem~\ref{th:1.7}
 will require a number of results about the structure of 
$H(\Gamma , \alpha)$ as an $\SS(\fg^*)$ module.  These results were proved
 in an earlier paper of ours on ``equivariant Morse theory on
 graphs'' (\cite{GZ3}), and we will refer to this paper for a detailed 
 treatment of the material in the next few paragraphs.  Let
 \begin{equation}
\label{eq:2.3}
   \P = \{ \xi \in \fg , \alpha_e (\xi) \neq 0
\hbox{ for all } e \in E_{\Gamma} \} \, .
 \end{equation}
Given an element, $\xi \in \P$, we will say that an oriented edge, 
$e$, is \emph{ascending with respect to} $\xi$ if $\alpha_e (\xi)>0$.
For every vertex, $p$, let $\sigma_p$, the \emph{index} of $p$,
be the number of ascending edges, $e$, with $t(e)=p$.

\begin{dfn}
  \label{dfn:2.5}
The $k^{th}$ Betti number, $b_{k} (\Gamma)$, is the number of
vertices, $p$, of $\Gamma$ for which $\sigma_p =k$.
\end{dfn}

\noindent \textbf{Remark: }
The definition of $\sigma_p$ depends upon the choice of $\xi$
but $b_{k}(\Gamma)$ turns out not to.  (See \cite[Theorem 2.6]{GZ1}).

\smallskip

A function $\phi : V_{\Gamma} \to \RR$ is a ($\xi$-compatible) 
\emph{Morse} function if, for every ascending edge $e$, 
$\phi (i(e))< \phi(t(e))$.  
It is not obvious, and in fact not true, that Morse
functions exist.  A necessary and sufficient condition for the
existence of a Morse function is that, for every ascending path
in $\Gamma$ the initial vertex of this path is distinct from its
terminal vertex (\emph{i.e.},~there are no ascending ``loops'').  If a
Morse function exists, however, one can easily perturb it so that 
it is \emph{injective} as a map of $V_{\Gamma}$ into $\RR$.  From 
now on we will let $\phi$ be a fixed Morse function with this
property.

The topological results discussed in Section~\ref{sec:1} prompt one to 
make the following Morse-theoretic conjectures about the
equivariant cohomology ring of a graph.

\begin{cnj}
\label{cnj:1}
$H(\Gamma ,\alpha)$ is a free $\SS(\fg^*)$ module with $b_{k} 
(\Gamma)$ generators of degree $k$.
\end{cnj}

\begin{cnj}
\label{cnj:2}
$H(\Gamma ,\alpha)$ is freely generated as an $S (\fg^*)$ 
module by a family of ``Thom classes''
$$\tau^+_p \in H^{k} (\Gamma ,\alpha), \,\,  k= \sigma_p \, ,   $$
satisfying
  \begin{equation}
    \label{eq:2.4}
    \hbox{support } \tau^+_p \subseteq   F_p
  \end{equation}
and
\begin{equation}
  \label{eq:2.5}
  \tau^+_p(p) = \prod_{e \in E^-_p} \alpha_e \,  \quad ( := \nu_p^+) \; ,  
\end{equation}
$F_p$ being the set of vertices which can be joined to $p$ by an
ascending path and $E^-_p$ being the set of descending edges in
$E_p$.
\end{cnj}

It is clear that Conjecture~\ref{cnj:2} implies Conjecture~\ref{cnj:1}, 
and it is not difficult to prove that Conjecture~\ref{cnj:1} implies 
Conjecture~\ref{cnj:2} (see \cite[\S2.4.3]{GZ2}). Therefore, since 
Conjecture~\ref{cnj:1} doesn't depend on the choice of an orientation 
of $\Gamma$ (\emph{i.e.} the choice of a polarizing vector $\xi \in \P$),
the same is true of Conjecture~\ref{cnj:2}. In particular, if we reverse 
the orientation  (replace $\xi$ by $-\xi$), we get from Conjecture~\ref{cnj:2}
the existence of Thom classes, $\tau_p^-$, $p \in V_{\Gamma}$, associated 
with the Morse function $-\phi$.

Unfortunately these conjectures are not true in general; 
however there is a useful necessary and 
 sufficient condition for them to be true involving
 certain subgraphs of $\Gamma$.

 \begin{dfn}
   \label{dfn:2.6}
A subgraph, $\Gamma_1$, of $\Gamma$ is \emph{totally geodesic}
if, for every pair of edges, $e$ and $e'$, of $\Gamma_1$, with
$i(e) =i(e')$, $\theta_e (e')$ is also an edge of $\Gamma_1$.
\end{dfn}

Note that if $\Gamma_1$ is a totally geodesic subgraph of
$\Gamma$ the restriction to it of $\alpha$ is, by \eqref{eq:2.2}, 
an axial function on $\Gamma_1$; so each of these subgraphs is
equipped with an action of $G$.  An important example of a totally 
geodesic subgraph is the following.  Let $\fh^*$ be a subspace of 
$\fg^*$, and let $\Gamma_{\fh^*}$ be the subgraph whose edges are 
the set
$$  \{ e \in E_{\Gamma}, \alpha_e \in \fh^* \} \, . $$
(It is clear, by \eqref{eq:2.1} and \eqref{eq:2.2} that this
\emph{is} totally geodesic.)  One of the main results of [GZ3] is
the following.

\begin{thm}
   \label{th:2.7}
Conjecture \ref{cnj:2} is true for $\Gamma$ if and only if, for every
two-dimensional subspace, $\fh^*$, of $\fg^*$, Conjecture~\ref{cnj:2}
 is true for $\Gamma_{\fh^*}$. 
\end{thm}

Thus, to verify that Conjecture~2 holds for $\Gamma$ it suffices
to verify it for these subgraphs (which is usually much easier
than verifying it for $\Gamma$ itself).

The proof of Theorem~\ref{th:2.7} involves a graph-theoretic
version of symplectic reduction.  We will say that $c$ is a
\emph{critical value} of the Morse function $\phi :V_{\Gamma} \to 
\RR$ if $c=\phi (p)$ for some $p \in V_{\Gamma}$ and, otherwise,
$c$ is a \emph{regular value}.  Let $c$ be a regular value of $\phi$ and 
let $V_c$ be the set of oriented edges, $e$, of $\Gamma$ with
$\phi (i(e))< c < \phi (t(e))$. We show in \cite{GZ3} that $V_c$ is the set 
of vertices of a hypergraph, $\Gamma_c$. Thus the elements of $V_c$ are 
both edges of the graph $\Gamma$ and vertices of this hypergraph. It is 
useful to distinguish between the two roles they play by saying that 
``an edge, $e$, intersects $\Gamma_c$ in a vertex, $v_e$.''

Let $\fg^*_{\xi}$ be the
annihilator of $\xi$ in $\fg^*$.  For each oriented edge, $e$, of 
$\Gamma$ we define a map $\rho_e : \fg^* \to \fg^*_{\xi}$ 
by setting
$$  \rho_e \alpha = \alpha -
  \frac{\alpha (\xi)}{\alpha_e (\xi)} \alpha_e \, .$$
This extends to a ring homomorphism
\begin{equation}
  \label{eq:2.6}
  \rho_e : \SS(\fg^*) \to \SS(\fg^*_{\xi})
\end{equation}
and, from \eqref{eq:2.6}, we get a ring homomorphism
$$  \K_c : H (\Gamma ,\alpha) \to \Maps(V_c , \SS(\fg^*_{\xi})) $$
by setting
\begin{equation}
\label{eq:2.7}
\K_c (g)(v_e) = \rho_e g_{i(e)} = \rho_eg_{t(e)} \, .
\end{equation}
(The two terms on the right are equal by \eqref{eq:1.10}.)

We show in \cite{GZ3} that $\K_c$ maps $H(\Gamma,\alpha)$ into the 
cohomology ring, $H(\Gamma_c , \alpha_c)$, of
the hypergraph $\Gamma_c$.  We won't bother to review here the
definition of this hypergraph cohomology ring (which is quite
tricky) since one of the main theorems of \cite{GZ3} asserts that, if
the hypotheses of Theorem~\ref{th:2.7} hold and if $\xi$
satisfies a certain genericity condition (which we will spell out 
below), the map
$$  \K_c : H(\Gamma ,\alpha) \to H (\Gamma_c ,\alpha_c)$$
is a submersion.  Hence, thanks to this theorem, one can
\emph{define} $H (\Gamma_c , \alpha_c)$ to be the image of
$\K_c$.

A key step in the proof of Theorem~\ref{th:2.7} is a theorem
which describes how the structure of the ring $H(\Gamma_c ,
\alpha_c)$ changes as one passes through a critical value of
$c$.  More explicitly suppose $c$ and $c'$ are regular values of
$\phi$ and suppose that there exists a unique vertex, $p$, with $
c<\phi (p)<c'$.  In addition suppose that, for $e_1$, $e_2$,
$e_3$, $e_4 \in E_p$
\begin{equation}
  \label{eq:2.8}
  \frac{1}{\alpha_{e_1}(\xi)} \rho_{e_2} \alpha_{e_1} \neq
  \frac{1}{\alpha_{e_3}(\xi)} \rho_{e_4} \alpha_{e_3}
\end{equation}
except when the two sides of \eqref{eq:2.8} are forced to be
equal (\emph{i.e.},~except when $e_1=e_2$ and $e_3=e_4$ or $e_1 =e_3$
and $e_2=e_4$).  The inequality \eqref{eq:2.8} is unfortunately
not satisfied for all elements, $\xi$, of the set \eqref{eq:2.3},
but one can show that those $\xi$'s for which it is satisfied
form an open dense subset of this set.

Let $r$ be the index of $p$ and let $s=d-r$.  Let $e_i, i=1,
\ldots ,r$ be the descending edges in $E_p$ and $e_a,
a=r+1 , \ldots ,d$ be the ascending edges in $E_p$.  Let
\begin{align}
  \Delta_c &= \{ e_i ; i=1, \ldots ,r\} \nonumber \\
\intertext{and}
\Delta_{c'} &=\{ e_a ; a=r+1,\ldots ,d \} \, .  \nonumber 
\end{align}
Then $\Delta_c$ is a subset of $V_c$, $\Delta'_c$ a subset of
$V_c$ and
$$  V_c-\Delta_c=V_{c'}-\Delta_{c'}=V_0 \; , $$
where $V_0$ is the intersection of $V_c$ and $V_{c'}$.  Let
$$  V^{\#} = V_0 \cup (\Delta_c \times \Delta_{c'}) \, .$$
Then one has projection maps
$$\pi_c  : V^{\#} \to V_c \quad \text{and} \quad \pi_{c'} :V^{\#} \to V_{c'}$$
and, from these projection maps, pull-back maps, $\pi^*_c$ and
$\pi^*_{c'}$, embedding the rings
\begin{equation}
  \label{eq:2.9}
  \Maps (V_c , \SS(\fg^*_{\xi}))
\end{equation}
and
\begin{equation}
  \label{eq:2.10}
  \Maps (V_{c'} , \SS(\fg^*_{\xi}))
\end{equation}
into the ring
\begin{equation}
  \label{eq:2.11}
  \Maps (V^{\#} , \SS(\fg^*_{\xi}))\, .  
\end{equation}
Moreover the ring
  \begin{equation}
    \label{eq:2.12}
    \Maps (\Delta_c , \SS(\fg^*_{\xi}))
  \end{equation}
sits in the ring \eqref{eq:2.9} as the set of maps $h: V_c \to 
\SS (\fg^*_{\xi})$  supported on $ \Delta_c$, and the ring
  \begin{equation}
    \label{eq:2.13}
    \Maps (\Delta'_c , \SS(\fg^*_{\xi}))
  \end{equation}
sits inside the ring \eqref{eq:2.10}; so \emph{all} the rings
\eqref{eq:2.9}--\eqref{eq:2.13} can be regarded as subrings of
\eqref{eq:2.11}.

Let $y_1 , \ldots , y_{n-1}$ be a basis of $\fg^*_{\xi}$ and $x$
a fixed element of $\fg^*$ with $\langle x,\xi \rangle =1$.  Let
\begin{align}
\alpha_{e_i} &=m_i (x-\beta_a(y)) \qquad i=1 ,\ldots , r \NR
\intertext{and}
\alpha_{e_a} &=m_a (x-\beta_a(y)) \qquad a=r+1 ,\ldots , d \nonumber \; ,  
\end{align}
with $m_i<0<m_a$ and with the $\beta$'s in $\fg^*_{\xi}$.
Consider the maps
\begin{align}
  \tau_c &: \Delta_c \to \fg^*_{\xi} \, ,& \, \tau_c (e_i) &=\beta_i \NR
\tau_{c'} &: \Delta_{c'} \to \fg^*_{\xi} \, ,& \, 
\tau_{c'} (e_a)&=\beta_a \NR
\intertext{and}
\tau^{\#} &: \Delta_c \times \Delta_{c'} \to \fg^*_{\xi} \, ,& \, 
\tau^{\#} (e_i,e_a)&= \beta_i- \beta_a \, . \nonumber
\end{align}
The first two of these maps depend on the choice of $x$; however, 
$\tau^{\#}$ is intrinsically defined since 
$\tau^{\#}(e_i,e_a)$ is just
$$  \frac{1}{\alpha_{e_i} (\xi)} \rho_{e_a} \alpha_{e_i} \, .$$
Also, by the genericity condition \eqref{eq:2.8} the map,
$\tau^{\#}$ sends $\Delta_c \times \Delta'_c$ injectively into
$\fg^*_{\xi}$ and, as a consequence, $\tau_c$ and $\tau'_c$ map $
\Delta_c$ and $\Delta_{c'}$ injectively into $\fg^*_{\xi}$.

Define the cohomology ring, $H (\Delta_c ,\tau_c)$, to be the set 
of all maps of $\Delta_c$ into $\SS(\fg^*_{\xi})$ of the form
$$  h=\sum^{r-1}_{i=0} g_i \tau^i_c \quad , \quad g_i \in \SS (\fg^*_{\xi}) $$
and define $H (\Delta_{c'}, \tau_{c'})$ to be the set
of all maps of $\Delta_{c'}$ into $\SS(\fg^*_{\xi})$ of the form
$$  h'=\sum^{s-1}_{i=0} g'_i \tau^i_{c'} \quad , \quad
      g'_i \in \SS (\fg^*_{\xi}) \, . $$
The theorem we alluded to above asserts

\begin{thm}
  \label{th:2.8}
For every $f \in H (\Gamma_c ,\alpha_c)$ and every $f_i \in H
(\Delta_c ,\tau_c)$, $i=1, \ldots ,s-1$, there exists a unique
$f' \in H (\Gamma_{c'}, \alpha_{c'})$ and  unique $f'_j \in H
(\Delta_{c'}, \tau_{c'})$, $j=1,\ldots , r-1$ such that
\begin{equation}
  \label{eq:2.14}
  f'+\sum^{r-1}_{j=1} (\tau^{\#})^j f'_j =
  f+\sum^{s-1}_{i=1} (\tau^{\#})^i f_i \, .
\end{equation}

\end{thm}

\noindent {\bf Remarks:}

\begin{enumerate}
\item 
This theorem gives one a concrete picture of how 
$H(\Gamma_c ,\alpha_c)$ changes as one goes through a critical point of
$\phi$.  Namely it shows that $H(\Gamma_{c'}, \alpha_{c'})$ can
be obtained from $H(\Gamma_c ,\alpha_c)$ by a ``blow-up''
followed by a ``blow-down''.  (Compare with \cite[Theorem~2.3.2]{GZ2}.)

\item 
This theorem can also be used to map cohomology classes in
$H(\Gamma_c ,\alpha_c)$ into cohomology classes in $H
(\Gamma_{c'}, \alpha_{c'})$.  By setting $f_i=0$, $i=1, \ldots
,s-1$, in \eqref{eq:2.14} one gets from the cohomology class
$f\in H(\Gamma_{c}, \alpha_{c})$ a unique cohomology class
$f'\in H(\Gamma_{c'}, \alpha_{c'})$.  (This observation will be
heavily exploited in the next section.)
\end{enumerate}

An important ingredient in the proof of Theorem~\ref{th:2.8} is the 
following.

\begin{thm}
\label{th:2.9}
If $f$ is in $H(\Gamma_c, \alpha_c)$, then its restriction to $\Delta_c$
is in $H(\Delta_c, \tau_c)$.
\end{thm}

(It is in the proof of this result that the hypotheses of 
Theorem~\ref{th:2.7} are needed.)

\section{Combinatorial formulas for Thom classes}\label{sec:3}

We will describe in this section how to compute the combinatorial 
Thom class, $\tau^+_{p_0}$, at an arbitrary point $p$ on the
flow-up $F_{p_o}$.  We recall that $\tau^+_{p_0}$ is canonically 
defined \emph{only} if the index function $\sigma : V_{\Gamma}
\to \ZZ$, is strictly increasing along ascending paths in
$\Gamma$.  Assuming that $\sigma$ has this property, we will show 
below that there is a simple inductive method for computing
$\tau^+_{p_0}$ on a critical level, $c$, of $\phi$ if one knows
the values of $\tau^+_{p_0}$ on lower critical levels.  Then,
later in this section, we will show that this method works even
when the hypothesis about $\sigma$ is dropped.  Let $\phi (p_0)=c_0$ and
$\sigma_{p_0}=m$.  The first step in this induction is to set
$\tau^+_{p_0}(q)=0$ for all vertices, $q$, with $\phi (q)<c_0$
and set $\tau^+_{p_0} (p_0)$ equal to $\nu_{p_0}^+$, as in 
\eqref{eq:2.5}.  
Now let $c>c_0$ and suppose, by induction, that $\tau^+_{p_0}(q)$ is
defined for all $q$ with $\phi (q) <c$ and is zero unless $q$ is
in $F_{p_0}$.  Let $p$ be a vertex with $\phi (p)=c$.  Let
$\sigma_p=r$ and let $e_{k}$, $k =1, \ldots , r$, be
the descending edges in $\Gamma$ with $p=i(e_{k})$.  Then
the vertices, $q_k=t(e_k)$, are points where
$\tau^+_{p_0}$ is already defined.  We will prove below.

\begin{lem}
\label{lem:3.1}
There exists a unique polynomial, $\psi \in \SS (\fg^*)$ such
that
\begin{equation}
  \label{eq:3.1}
  \psi \equiv \tau_{p_0}(q_{k}) \mod \alpha_{e_{k}}, \quad
  k =1 , \ldots , r \, .
\end{equation}
\end{lem}

\noindent{\textbf{Remark: }} 
The ``uniqueness'' part of this lemma is where the hypothesis on
$\sigma$ is used. If $f \in \SS^m (\fg^*)$ and 
\begin{align}
  f &= 0 \mod \alpha_{e_{k}}, \quad  k =1, \ldots , r \; , \NR
\intertext{then}
  f &= h \alpha_{e_1} \ldots \alpha_{e_r}, \quad h \in 
\SS^{m-r}(\fg^*)\, . \nonumber
\end{align}
Hence, if $m<r$, $f$ is identically zero.

Using this lemma, set $\tau^+_{p_0}(p) =\psi$, and continue
with the induction until the set of vertices of $\Gamma$ is
exhausted.  It is clear from \eqref{eq:3.1} that this
construction gives us a map:  
$\tau^+_{p_0}:V_{\Gamma} \to \SS^m(\fg^*)$ 
satisfying \eqref{eq:1.10} and that this map is supported on $F_p$.

By giving a constructive proof of the ``existence'' part of
Lemma~\ref{lem:3.1} the induction argument we just sketched can be
converted into a formula for $\tau^+_{p_0}$, and this will be the 
main goal of this sections.  Note that the solution
of \eqref{eq:3.1} is basically an interpolation problem:  finding 
a polynomial with prescribed values at 
$\alpha_{e_1}, \ldots ,\alpha_{e_r}$.  
To solve this problem constructively, we review
a few elementary facts about ``interpolation''.  The basic
problem in interpolation theory is to find a polynomial
\begin{equation}
  \label{eq:3.2}
  \fp (x) = \sum^n_{i=1} g_i x^{i-1}
\end{equation}
which takes prescribed values
\begin{equation}
  \label{eq:3.3}
  \fp (x_i)=f_i
\end{equation}
at $n$ distinct points, $x_i$, on the complex line.  The solution 
of this problem is more or less trivial.  The polynomial
\begin{equation}
  \label{eq:3.4}
  \fp (x) = \sum_j \prod_{k \neq j}
  \frac{x-x_{k}}{x_j-x_{k}} f_j
\end{equation}
satisfies \eqref{eq:3.3} and is the only polynomial of
degree less than $n$ which \emph{does} satisfy \eqref{eq:3.3}.  Moreover, 
from \eqref{eq:3.4} one gets an explicit formula for the $g_i$'s
in \eqref{eq:3.2}.  Let
$$\prod_{\ell \neq j} (x-x_{\ell}) = \sum^{n-1}_{i=1}
  (-1)^{n-1-i}\sigma^j_{n-1-i} x^i \; , $$
where $\sigma^j_r$ is the $r$-th elementary symmetric function
in the variables, $x_{\ell}$, $\ell \neq j$.  Then
by \eqref{eq:3.4}
\begin{equation}
  \label{eq:3.5}
  g_i = (-1)^{n-i} \sum_{j=1}^n 
\frac{\sigma^j_{n-i}}{\prod_{\ell \neq j} (x_j-x_{\ell})} f_j \, .
\end{equation}
One consequence of \eqref{eq:3.5} is an inversion formula for the 
Vandermonde matrix, $A$, with entries
$$  a_{ij} = x^{j-1}_i \, , \,\,
  1 \leq j,i \leq n \, .$$
If $B=A^{-1}$ then by \eqref{eq:3.5} and \eqref{eq:3.3}:
\begin{equation}
  \label{eq:3.6}
  b_{ij} = (-1)^{n-i} 
\frac{\sigma^j_{n-i}}{\prod_{\ell \neq j} (x_j-x_{\ell})} \, .
\end{equation}
In particular
\begin{align}
\label{eq:3.7}
b_{nj} &=  \prod_{\ell \neq j} \frac{1}{x_j - x_{\ell}}\\
\intertext{and}
\label{eq:3.8}
b_{1j} &= \prod_{\ell \neq j} \frac{-x_{\ell}}{x_j-x_{\ell}} \, .
\end{align}
It is sometimes convenient to write the inversion
formula \eqref{eq:3.6} in terms of the elementary symmetric
functions  $\sigma_r = \sigma_r (x_1 , \ldots , x_n)$ rather
than in terms of the $\sigma^j_r$'s.  To do so, we note that
\begin{equation}
  \label{eq:3.9}
  \sigma^j_{k} = \sum^{k}_{r=0} (-1)^r  \sigma_{k -r} x^r_j \, .
\end{equation}
(To derive \eqref{eq:3.9} observe that 
$$\prod_{\ell \neq j} (x-x_{\ell}) = 
\prod_{\ell} (x-x_{\ell}) \frac{1}{x-x_j} 
= \prod_{\ell} (x-x_{\ell}) \frac{1}{x}
\sum^{\infty}_{i=0} \Bigl( \frac{x_j}{x} \Bigr)^i$$
and compare coefficients of $x^{n-k-1}$ on both sides.)
Substituting \eqref{eq:3.9} into \eqref{eq:3.6} one gets an
alternative inversion formula for the Vandermonde matrix
\begin{equation}
  \label{eq:3.10}
  b_{ij} = \sum^{n-i}_{r=0} (-1)^{n-i-r} 
\frac{\sigma_{n-i-r} x^r_j}{\prod_{\ell \neq j} (x_j-x_{\ell})} \, .
\end{equation}
 Finally we note a couple of trivial consequences of \eqref{eq:3.7}
 and \eqref{eq:3.8}.  From \eqref{eq:3.7} and the identity
$$ \sum_j b_{nj} a_{jk} = \delta^n_{k}$$
we conclude that the sum
$$  \sum_j \frac{x_j^{k -1}}
      {\prod_{\ell \neq j}(x_j-x_{\ell})}$$
is zero if $k$ is less than $n$ and $1$ if $k =n$; and
from \eqref{eq:3.8} and the identity
$$  \sum_j b_{1j} a_{jk} =1$$
we conclude that
\begin{equation}
  \label{eq:3.11}
  \sum_j \prod_{\ell \neq j} \frac{-x_{\ell}}{x_j-x_{\ell}} =1 \, .
\end{equation}

In the applications which we will make of these identities below
the $x_i$'s will be indeterminants and the $f_i$'s polynomials in these
indeterminants, and we will want to know when the
$g_i$'s are also polynomials in these indeterminants.  To answer
this question we will show that these identities have a simple
``topological'' interpretation:  Suppose one is given a graph,
$\Gamma$, and an action of $G$ on $\Gamma$ defined by an axial
function, $\alpha: E_{\Gamma} \to \fg^*$.  One of the main
results of an earlier paper of ours is that there is a canonical
integration operation
$$  \int_{\Gamma} : H (\Gamma ,\alpha) \to S(\fg^*)$$
defined by
\begin{equation}
  \label{eq:3.12}
  \int_{\Gamma} f = \sum_{p \in V_{\Gamma}} f_p \delta_p
\end{equation}
where
$$\delta_p =\bigl( \prod_{i(e) =p} \alpha_e \bigr)^{-1}\, .$$
(See \cite[\S~2.4]{GZ1}. This formula is the formal analogue of the
standard localization theorem \cite{AB}-\cite{BV} in equivariant DeRham 
theory.)   

In particular let $\Delta$ be the complete
graph on $n$ vertices.  Denote these vertices by $1,\ldots ,n$,
and let $x_1 ,\ldots ,x_n$ be a basis of $\fg^*$.  It is easy to
check that the map 
$$  \alpha : E_{\Delta} \to \fg^* \, , $$
which assigns the weight $x_i-x_j$ to the edge joining $i$ to
$j$, is an axial function, and that the map
$$\tau : V_{\Delta} \to \fg^* , \,\, \tau(i) = x_i$$
is an element of $H^1 (\Delta , \alpha)$.  We claim that 
$1,\tau, \ldots , \tau^{n-1}$ generate $H(\Delta, \alpha)$ as a module
over $\SS (\fg^*)$.  To see this let $\nu_i$ be the cohomology
class
$$  \nu_i = \sum_{r=0}^{n-i} (-1)^{n-i-r} \sigma_{n-i-r} \tau^r \, ,\,\,
  i=1, \ldots ,n \, .$$
Then \eqref{eq:3.10} simply asserts that
$$  \int_{\Delta} \nu_i \tau^{j-1} = \delta^j_i \, .$$
In particular if $f$ is any cohomology class, then one can 
express $f$ as a sum
$$  f = \sum_{i=1}^{n} g_i \tau^{i-1}$$
where
$$g_i = \int_{\Delta} \nu_i f \, \in \SS (\fg^*) \, .$$

This proves the assertion:

\begin{prop}
  \label{prop:3.2}
  If $f_1 , \ldots , f_n$ are polynomials in $x_1 , \ldots , x_n$ 
  and the function
$$\fp (x) = \sum^n_{i=1} g_i x^{i-1}$$
solves the interpolation problem
$$\fp (x_i) = f_i$$
then the $g_i$'s are polynomials in $x_1 , \ldots , x_n$ if and
only if $x_i-x_j$ divides $f_i-f_j$.
\end{prop}

\begin{proof}
  If $x_i-x_j$ divides $f_i-f_j$ the map
$$f: V_{\Delta} \to S (\fg^*), \,\, i \to f_i$$
is in $H (\Delta , \alpha)$.
\end{proof}

Let's come back now to Theorem~\ref{th:2.8} and the application
of it which we discussed at the end of Section~\ref{sec:2}.  As in
Theorem~\ref{th:2.8} let $c$ and $c'$ be regular values of
$\phi$, and suppose that there is just one vertex, $p$, of
$\Gamma$ with $c< \phi (p) < c'$.  By setting 
$f_1 = f_2 = \cdots f_{s-1}=0$ in \eqref{eq:2.14} one gets a map
\begin{equation}
  \label{eq:3.13}
  T_{c,c'} : H (\Gamma_c , \alpha_c) \to H (\Gamma_{c'},\alpha_{c'})
\end{equation}
sending $f_0$ to $f'_0$, and by the results above one can give a
fairly concrete description of this map.  Let's order the edges $
e_1, \ldots ,e_d \in E_p$ so that $e_1 , \ldots ,e_r$ are
descending and $e_{r+1}, \ldots , e_d$ are ascending, and let
$\Delta_c$ and $\Delta'_c$ be the vertices of $\Gamma_c$ and
$\Gamma'_c$ corresponding to the $e_j$'s, $1 \leq j \leq r$, and
the $e_a$'s, $r+1 \leq a \leq d$.  Then
\begin{align}
  V_c &= V_0 \cup \Delta_c  \nonumber \\
\intertext{and}
V_{c'} &= V_0 \cup \Delta_{c'}  \; , \nonumber
\end{align}
$V_0$ being the vertices which are common to $\Gamma_c$ and
$\Gamma_{c'}$.  To simplify notation we will identify $\Delta_c$
with the set $\{ 1,\ldots ,r \}$ and $\Delta_{c'}$ with the set 
$\{ r+1 , \ldots , d \}$.  Let $f_0$ be in 
$H (\Gamma_c,\alpha_c)$ and let $f'_0 = T_{c,c'} (f_0)$.  Then, by
\eqref{eq:2.14} and \eqref{eq:3.8}
\begin{equation}
  \label{eq:3.14}
  f'_0 (a) = \sum^r_{j=1} \prod_{k \neq j}
  \frac{\beta_a - \beta_{k}}{\beta_j - \beta_{k}}
  f_0 (j) \, ,
\end{equation}
for $a$ in $\Delta_{c'}$ and $j$ and $k$ in $\Delta_c$; and
\begin{equation}
  \label{eq:3.15}
  f'_0 = f_0 \quad \hbox{ on } V_0 \, .
\end{equation}

The identity \eqref{eq:3.14} has the following simple
interpretation.  Let
\begin{equation}
  \label{eq:3.16}
  \fp (x) = \sum^r_{j=1} \prod_{k \neq j}
  \frac{x-\beta_{k}}{\beta_j - \beta_{k}} f_0 (j) \, .
\end{equation}
Then, by \eqref{eq:3.4}, $\fp (x)$ solves the interpolation problem
\begin{equation}
  \label{eq:3.17}
  \fp (\beta_j) = f_0 (j) \, .
\end{equation}
On the other hand, by Theorem~\ref{th:2.9}
$$f_0 |_{ \Delta_c} \in H(\Delta_c , \tau_c) \, ;$$
so, by Proposition~\ref{prop:3.2}, $\fp(x)$
 is a polynomial in $x, \beta_1 (y), \ldots , \beta_r (y)  $ and
 hence also a polynomial in $(x,y_1 , \ldots , y_{n-1})$,
 \emph{i.e.},~an element of the ring, $\SS (\fg^*)$.  In fact, if 
$f_0 \in  H^m (\Gamma_c , \alpha_c)$ and $r>m$, $\fp (x)$ is the
 \emph{unique} element of $\SS^m (\fg^*)$ satisfying
 \eqref{eq:3.17}.  Now by \eqref{eq:3.14}, $\fp (\beta_a)= f'_0
 (a)$, so \eqref{eq:3.14} simply says that
$$f'_0 |_{ \Delta_{c'}} = \fp (\tau_{c'}) \, .$$
Thus to summarize, we have proved:

\begin{thm}
  \label{th:3.3}
The map
$$T_{c,c'} : H^m (\Gamma_c , \alpha_c) \to H^m (\Gamma_{c'}, \alpha_{c'})$$
is the identity map on $V_0$ and on $\Delta_c$ is the
``flip--flop''
\begin{equation}
  \label{eq:3.18}
  f_0 =\fp (\tau_c) \to \fp (x) \to f'_0 =\fp (\tau_{c'}).
\end{equation}

\end{thm}

Since $V_c$ and $V_{c'}$ are finite sets, the map $T_{c,c'}$ is 
defined by a matrix with entries
$$  T_{c,c'} (v,v'), \quad (v,v') \in V_c \times V_{c'} \, .$$
An important property of this matrix is the Markov property:
\begin{equation}
  \label{eq:3.19}
  \sum_{v \in V_c} T_{c,c'} (v,v') =1 \, .
\end{equation}

\begin{proof}
  It suffices to check this for $a \in \Delta_{c'}$, \emph{i.e.} it 
suffices to check that
$$  \sum^r_{j=1} T_{c,c'} (j,a) =1 \, .$$
However, by \eqref{eq:3.14}, this sum is equal to
$$  \sum^r_{j=1} \prod_{k \neq j} 
  \frac{\beta_a - \beta_{k}}{\beta_j - \beta_{k}}$$
which is equal to $1$ by \eqref{eq:3.11}, 
with $x_{\ell} = \beta_{\ell} -\beta_a$.
\end{proof}

We will next give a more intrinsic description of $T_{c,c'}$ 
and of the polynomial $\fp$ in \eqref{eq:3.16}.  We recall that
\begin{align}
  \alpha_{e_i} &= m_i (x-\beta_i (y)), \quad i=1, \ldots ,r \nonumber \\
\intertext{and}
\alpha_{e_a} &= m_a (x-\beta_a(y)), \quad a=r+1, \ldots ,d \, .  \nonumber
\end{align}
Hence
$$T_{c,c'} (j,a) = \prod_{k \neq j}
\frac{\beta_a - \beta_{k}}{\beta_j - \beta_{k}} =
 \prod_{k \neq j}
\frac{\alpha_{e_{k}}-(m_{k} /m_a) \alpha_{e_a}}
              {\alpha_{e_{k}} - (m_{k}/m_j) \alpha_{e_j}}$$
and therefore
$$  T_{c,c'} (j,a) = 
  \frac{\rho_{e_a} (\prod_{k \neq j}\alpha_{e_{k}})}
       {\rho_{e_j} (\prod_{k \neq j} \alpha_{e_{k}})} \; , $$
where $\rho_e$ is the map \eqref{eq:2.6}.  Similarly the polynomial
$\fp$ is just
\begin{equation}
  \label{eq:3.20}
  \sum_j \frac{\prod_{k \neq j} \alpha_{e_{k}}}
     {\rho_{e_j} (\prod_{k \neq j} \alpha_{e_{k}})} f_0 (j) \, .
\end{equation}

By iterating \eqref{eq:3.13} we will extend the definition of
$T_{c,c'}$ to arbitrary regular values of $\phi$ with $c <c'$.
Let
$$ c_i , i=0 , \ldots , \ell $$
be regular values of $\phi$ with $c_0 = c$ and $c_{\ell} =c'$
such that there exists a unique vertex, $p_i$, with 
$c_{i-1} < \phi (p_i) < c_i$.  Let $T_i = T_{c_{i-1}}, c_i$ and let
\begin{equation}
  \label{eq:3.21}
  T : H (\Gamma_c ,\alpha_c) \to H(\Gamma_{c'},\alpha_{c'})
\end{equation}
be the map
\begin{equation}
  \label{eq:3.22}
  T= T_{\ell} \circ \cdots \circ T_1 \, .
\end{equation}
We will list a few properties of this map:

\bigskip

{\bf 1). } 
This map is defined by a matrix with entries
$$T (v,v') , \,\, (v,v') \in V_c \times V_{c'}$$
and since all the factors on the right hand side of
\eqref{eq:3.22} have the Markov property \eqref{eq:3.19}, this
matrix also has this property.

\bigskip

{\bf 2). } 
  The matrix version of \eqref{eq:3.22} asserts that 
$$T (v,w) = \sum T_{\ell}(v_{\ell -1},w) \cdots T_2 (v_1,v_2) T_1 (v,v_1)$$
summed over all sequences $v_1 , \ldots , v_{\ell -1}$ with
$v_{k} \in V_{c_{k}}$.  By \eqref{eq:3.14} and \eqref{eq:3.15}, a large
number of the matrix entries in this formula are either $1$ or
$0$:  If $e_{k} $ is an ascending edge which intersects
$\Gamma_{c_{k}}$ in $v_{k}$ and $\Gamma_{c_{k-1}}$ in $v_{k-1}$, then 
$$T(v', v_k) = \begin{cases}
0 & \text{  if  } v' \neq v_{k-1} , \\ 
1 & \text{ if } v' = v_{k-1}  \end{cases} \; . $$
This fact can be exploited to
write the sum above more succinctly.  For every pair of edges, 
$e$ and $e'$, with $t(e) =i(e') =p$, let $e_1 , \ldots ,e_r$
be the descending edges in $E_p$, ordered so that $e_r =\bar{e}$
and let 
$$  Q (e,e') = \frac{\rho_{e'} (\prod^{r-1}_{i=1} \alpha_{e_i})}
     {\rho_e (\prod^{r-1}_{i=1} \alpha_{e_i}) } \, . $$
Then $T(v,w) $ can be written as a weighted sum:
$$  T(v,w) = \sum_{\gamma} Q(\gamma) $$
over all ascending paths, $\gamma$, in $\Gamma$ whose initial
edge intersects $\Gamma_c$ in $v$ and whose terminal edge
intersects $\Gamma_{c'}$ in $w$, the weighting of the path, $\gamma$,
being given by
\begin{equation}
  \label{eq:3.23}
  Q (\gamma) = \prod^m_{i=1} Q(e_{i-1},e_i) \, ,
\end{equation}
where the $e_i$'s are the edges of $\gamma$, ordered so that 
for $i >1$, $t(e_{i-1}) =i(e_i)$.

\bigskip

{\bf 3). } 
  The map \eqref{eq:3.22} can also be viewed as a series of
  ``flip--flops''.  Let $f_0$ be an element of $H (\Gamma_c ,
  \alpha_c) $ and let $f_i=T_i \ldots T_1 f$.  Then $T_i$ maps
  $f_{i-1}$ to $f_i$ by a map of the form \eqref{eq:3.18}.  Let's 
  denote the polynomial, $\fp$, in \eqref{eq:3.18} by
  $\psi_{p_i}$.  We claim:

  \begin{prop}
    \label{prop:3.4}
If $p_i$ is joined to $p_j$ by an ascending edge, $e$,
$$  \psi_{p_i} \equiv \psi_{p_j} \mod \alpha_e \, . $$
\end{prop}

\begin{proof}
    This is equivalent to asserting that
    \begin{equation}
      \label{eq:3.24}
      \rho_e \psi_{p_i} = \rho_e \psi_{p_j} \, ;
    \end{equation}
however, \eqref{eq:3.24} is, by definition, the common value of $
f_{k} (v_{k})$, $i \leq k <j$, at the vertices,
$v_{k}$, at which $e$ intersects $\Gamma_{c_{k}}$.

\bigskip

{\bf 4). } 
  In particular let $p_0$ be an arbitrary vertex of $\Gamma$; and 
  choose $c$ and $c'$ such that there are no critical values of
  $\phi$ on the interval, $(\phi (p_0),c)$ and such that $c' >
  \max \phi (p)$, $p \in V_{\Gamma}$. 
Order the edges $e_1 , \ldots ,e_d$ in $E_{p_0}$ so that 
$e_1 , \ldots ,e_r$ are 
descending and $e_{r+1}, \ldots ,e_d$ are ascending.  For $r+1
  \leq a \leq d$ let $v_a$ be the vertex at which $e_a$
  intersects $\Gamma_c$ and let
$$ f_0 : V_{\Gamma_c} \to \SS^r (\fg^*_{\xi})$$
be the map defined by
\begin{equation}
\label{eq:3.25}
f_0 (v) = \begin{cases} 
0 , & \text{ if } v \not\in \{ v_{r+1}, \ldots , v_d \} \\
\rho_{e_a} (\prod^r_{i=1} \alpha_{e_i}),  & \text{ if } v = v_a.
\end{cases}
\end{equation} \qed
\renewcommand{\qed}{}
\end{proof}

\begin{prop}
\label{prop:3.5}
$f_0$ is an element of $H(\Gamma_c ,\alpha_c)$.
  \end{prop}

\begin{proof}
    By \eqref{eq:2.7}, $f_0 = \K_c \tau^+_{p_0}$.  (This
    proof assumes that there \emph{exists} a Thom class,
    $\tau^+_{p_0}$, having the properties listed in
    Theorem~\ref{th:1.5}.  Alternatively, 
Proposition~\ref{prop:3.5} can be proved directly using a
    more sophisticated definition of $H (\Gamma_c ,\alpha_c)$
    than that which we gave in Section~\ref{sec:2}.  For more details
    see \cite[\S~4]{GZ3}.)
\end{proof}

By applying the sequence of flip--flops, $T_i$, to the $f_0$
above, we get a polynomial, $\psi_{p_i} \in \SS^r (\fg^*)$ for each 
vertex, $p_i$, of $\Gamma$ with $\phi (p_i) >c$.  On the other
hand, we can define $\tau_p$ for $\phi (p) <c$ to be equal to
\eqref{eq:2.5} at $p_0$ and equal to zero otherwise.  By
\eqref{eq:3.24} $\tau_p$ satisfies the cocycle condition
\eqref{eq:1.10} at all vertices except $p_0$, and by
\eqref{eq:3.25} it satisfies this condition at $p_0$ as well.
Thus, if the index function, $\sigma : V_{\Gamma} \to \ZZ$, is
strictly increasing along ascending paths, this settles the
existence part of Lemma~\ref{lem:3.1} and justifies the
induction method for constructing $\tau^+_{p_0} $ which we
outlined at the beginning of this section.  On the other hand, if 
$\sigma$ fails to satisfy this hypothesis, the assignment, $p \to 
\tau_p$, still defines an element of $H^r (\Gamma ,\alpha)$ with
the properties listed in Theorem~\ref{th:1.5}; however, it
won't be the only element with these properties and may not even 
be the optimal element with these properties.

\bigskip

{\bf 5). } 
  From \eqref{eq:3.23} one gets the following ``path integral''
  formula for $\tau^+_p$.  If $e$ is an ascending edge of
  $\Gamma$, let $p=t(e)$ and let $e_1 , \ldots , e_r$ be the
  descending edges in $E_p$, ordered so that $e_r =\bar{e}$.  Let
  \begin{equation}
    \label{eq:3.26}
    Q(e) = \frac{\prod^{r-1}_{i=1} \alpha_{e_i}}
      {\rho_e (\prod^{r-1}_{i-1} \alpha_{e_i})} \, .
  \end{equation}
Then by \eqref{eq:3.20}, \eqref{eq:3.23} and \eqref{eq:3.25}
\begin{equation}
  \label{eq:3.27}
  \tau^+_{p_0} (p) = \sum E (\gamma)
\end{equation}
summed over all ascending paths in $\Gamma$ joining $p_0$ to $p$, 
$E(\gamma)$ being defined by
\begin{equation}
  \label{eq:3.28}
  E (\gamma) = Q (e_m) Q(\gamma) \rho_{e_1}(\nu_{p_0}^+) \; , 
\end{equation}
where $e_1$ is the initial edge of $\gamma$ and $e_m$ is the
terminal edge of $\gamma$.

\section{Combinatorial intersection numbers}\label{sec:4}

We will show below how to recast the formula \eqref{eq:3.28} into the form
\eqref{eq:1.12} and will also show that, if the hypothesis of 
Theorem~\ref{th:1.6} is satisfied, one can deduce from \eqref{eq:1.12} 
the formula that we described in Section 1 for the products of Thom classes. 
First, however, we will examine this hypothesis in more detail: Suppose the 
graph $\Gamma$ is connected and admits a family of Thom classes, 
$\tau_p^+$, $p \in V_{\Gamma}$, which generates $H(\Gamma, \alpha)$ 
as a free module over the ring $\SS(\fg^*)$, and have the properties 
\eqref{eq:2.4} and \eqref{eq:2.5}. By \eqref{eq:1.10} 
$$\dim H^0(\Gamma, \alpha) = 1 \; ; $$
hence there is a unique vertex, $p_0$, with $\sigma_{p_0} = 0$. Let $p$ be an 
arbitrary vertex of $\Gamma$ and let $\gamma$ be an ascending path with 
terminal endpoint $p$. If $\gamma$ is of \emph{maximal} length, its initial 
vertex has to be $p_0$, since every other vertex has a descending edge. 
Let $\phi(p)$ be the length of this longest path. If $p$ can be joined to $q$
by an ascending edge, $\phi(p)$ is strictly less than $\phi(q)$, so the map
$$\phi : V_{\Gamma} \to \ZZ \quad , \quad p \to \phi(p) \; , $$
is a Morse function.

\begin{thm}\label{th:4.1}
The index function, $\sigma$, is strictly increasing along ascending paths 
if and only if $\phi = \sigma$ (\emph{i.e.} if and only if the Morse 
function, $\phi$, is self-indexing.) Moreover, if $\phi$ has this property 
then, for every pair of vertices, $p \in V_{\Gamma}$ and $q \in F_p$, the 
length of the longest ascending path from $p$ to $q$ 
is $\sigma_q-\sigma_p$.
\end{thm}

It suffices to prove the last assertion, and it suffices by induction to 
prove this assertion for paths of length one. This we will do by proving a 
slightly stronger assertion.

\begin{thm}
\label{th:4.2}
Let $e$ be an ascending edge joining $p$ to $q$. If 
$e$ is the only ascending path from $p$ to $q$ then 
$$\sigma_q \leq \sigma_p+1 \; .$$
\end{thm}

\begin{proof}
Let $\Gamma_e$ be the totally geodesic subgraph of $\Gamma$ consisting 
of the single edge, $e$, and vertices $p$ and $q$. The Thom class, $\tau_e$, 
of $\Gamma_e$ is defined by
$$\tau_e(p) = \prod_{\substack{i(e')=p \\ e' \neq e}} \alpha_{e'} 
\quad\; , \qquad
\tau_e(q) = \prod_{\substack{i(e'')=q \\ e'' \neq \bar{e}}} \alpha_{e''} $$
and 
$$\tau_e(r) =0 \qquad \text{ if } r \neq p,q \; . $$
It is easily checked that $\tau_e \in H^{d-1}(\Gamma, \alpha)$.

\begin{lem}
A cohomology class, $\tau \in H(\Gamma, \alpha)$ is supported on $\{ p,q \}$ 
iff $\tau = h \tau_e$, $h \in H(\Gamma_e, \alpha)$.
\end{lem}

Suppose now that $e$ satisfies the hypotheses of Theorem~\ref{th:4.2}. Then
$\tau_p^+ \tau_q^-$ is supported on $\{ p,q \}$; so, by the lemma, 
\begin{equation}
\label{eq:4.1}
\tau_p^+ \tau_q^- =  h \tau_e \; , \; h \in H(\Gamma_e, \alpha) \; .
\end{equation}
In particular
$$\sigma_p + d -\sigma_q = \text{degree }\tau_p^+ + \text{degree }\tau_q^- 
\geq \text{degree } \tau_e  =  d-1 \; , $$
so $\sigma_q \leq \sigma_p +1$.
\end{proof}

Coming back to the formula \eqref{eq:3.27} lets first consider the simplest 
summands in this formula, those associated with paths, $\gamma$, of 
length one. For each $q \in V_{\Gamma}$ denote by $E_q^-$ and $E_q^+$ 
the descending and ascending edges in $E_q$ and let $\nu_q$ be defined as 
in \eqref{eq:2.5}. Let $\gamma$ be an ascending path of length one 
consisting of a single edge, $e$, with $i(e)=p $ and $t(e) =q$. 
Then by \eqref{eq:3.26} and \eqref{eq:3.28}
\begin{equation}
\label{eq:4.2}
E(\gamma) = \frac{\nu_q}{-\alpha_e} \cdot 
\frac{\prod{}{'} \rho_e(\alpha_{e_i})}
{\prod{}{''} \rho_e(\alpha_{e_j'})} \; ,
\end{equation}
where $\prod'$ in the enumerator is a product 
over the edges $e_i \in E_p^-$ and $\prod''$ in the denominator 
is the product over the edges $e_j' \in E_q^- -\{ qp \}$. 
Let 
$$\theta_e : E_p \to E_q$$
be the connection along this edge and let
$\theta_{\bar{e}} = \theta_e^{-1} : E_q \to E_p$. We define
\begin{align}
\label{eq:4.3}
E_{p,q} & = \{ e' \in E_p^- \; ; \; \theta_e(e') \not\in E_q^- \} \\
\intertext{and}
\label{eq:4.4}
E_{q,p} & = \{ e'' \in E_q^- \; ; \; \theta_{\bar{e}}(e'') 
\not\in E_p^- \} - \{ \bar{e} \} \; .  
\end{align}
Note that $\theta_e$ restricts to a bijection 
$$\theta_e : E_p - E_{p,q}  \to E_q - E_{p,q} \; .$$ 
If $e' \in E_p$, then \eqref{eq:2.2} implies 
\begin{equation}
\label{eq:4.5}
\rho_e(\alpha_{e'}) = \rho_e(\alpha_{\theta_e(e')}) \; .
\end{equation}
Therefore if $e_i \in E_p -E_{p,q}$, then the terms 
corresponding to $e_i$ and $\theta_e(e_i)$ in \eqref{eq:4.2} 
cancel each other and we obtain
$$\frac{E(\gamma)}{\nu_p} = \frac{\nu_q}{-\alpha_e \nu_p}  \cdot 
\frac{\rho_e(Z_{p,q})}{\rho_e(Z_{q,p})} = 
\frac{- \nu_q}{\alpha_{pq} \nu_p} \cdot \Theta_{pq} \; ,$$
where 
$$Z_{p,q} = \prod_{e' \in E_{p,q}} \alpha_{e'} \quad , \quad 
Z_{q,p} = \prod_{e'' \in E_{q,p}} \alpha_{e''} \; , $$
and 
\begin{equation}
\label{eq:4.6}
\Theta_{pq} = 
\frac{\prod{}{'} \rho_e(\alpha_{e_i})}
{\prod{}{''} \rho_e(\alpha_{e_j'})} =
\frac{\rho_e(Z_{p,q})}{\rho_e(Z_{q,p})} \; . 
\end{equation}

If $\gamma$ is the only ascending path from $p$ to $q$, 
then $\Theta_{p,q}$ has an interpretation as an 
``intersection number'': By \eqref{eq:4.1} the quotient,
$$\frac{\tau_p^+ \tau_q^-}{\tau_e}$$
is an element of $H(\Gamma_e, \alpha_e)$. Let $c$ be a point on the interval
$(\phi(p) , \phi(q))$ and let $v_e$ be the vertex of $\Gamma_c$ corresponding
to $e$. If we apply the Kirwan map
$$\K_c : H(\Gamma, \alpha) \to H(\Gamma_c, \alpha_c)$$
to this quotient and evaluate at $v_e$ we get an element of 
$\SS(\fg_{\xi}^*)$.
We claim
\begin{equation}
\label{eq:4.7}
\frac{\Theta_{p,q}}{\alpha_e(\xi)} = 
\K_c \Bigl( \frac{\tau_p^+ \tau_q^-}{\tau_e} \Bigr) (v_e) \; . 
\end{equation}

\begin{proof}
A direct computation shows that 
$$\K_c (\tau_p^+) (v_e) = \prod{}{'} \rho_e(\alpha_{e'}) \quad
\text{ and } \quad 
\K_c(\tau_q^-)(v_e)= \frac{\K_c (\tau_e) (v_e)}
{\prod'' \rho_e (\alpha_e'')} \; , $$
hence, \eqref{eq:4.7} follows from \eqref{eq:4.6}. 
\end{proof}

We will now show that the right hand side of \eqref{eq:4.7} can be 
interpreted as a ``pairing'' of the cohomology classes $\K_c(\tau_p^+)$ 
and $\K_c(\tau_q^-)$.  We pointed out in Section~\ref{sec:3} that the 
localization formula in equivariant DeRham theory enables one to define 
an integration operation on $H(\Gamma, \alpha)$. The analogue of this 
result for $\Gamma_c$ asserts that there is an integration operation
$$\int_{\Gamma_c} : H(\Gamma_c, \alpha_c ) \to \SS(\fg_{\xi}^*)$$
mapping $f \in H(\Gamma_c, \alpha_c)$ to the sum
$$\sum_{v \in V_c} f(v) \delta_v \; , $$
where
$$\delta_v = \bigl( \K_c(\tau_e)(v) \bigr)^{-1} \; , $$
$e$ being the edge of $\Gamma$ which intersects $\Gamma_c$ of the vertex 
$v = v_e$.

In particular, consider the product in $H(\Gamma_c, \alpha_c)$ of 
$\K_c(\tau_p^+)$ and $\K_c(\tau_q^-)$. If $e$ is the only ascending 
path in $\Gamma$ joining $p$ to $q$ this product is zero except at the point 
$v_e$; so by \eqref{eq:4.7}
\begin{equation}
\label{eq:4.8}
\frac{\Theta_{p,q}}{\alpha_e(\xi)} = 
\int_{\Gamma_c} \K_c(\tau_p^+)\K_c(\tau_q^-) \; , 
\end{equation}
which is the formal analogue of the intersection number \eqref{eq:1.4}.

\smallskip

\noindent
{\bf Remarks:}

\begin{enumerate} 

\item  By Theorem~\ref{th:4.2}, $\sigma_q \leq \sigma_p +1$. One can see by 
inspection that the right hand side of \eqref{eq:4.8}, which is by 
definition an element of $\SS(\fg_{\xi}^*)$, is of degree 
$\sigma_p+1-\sigma_q$. In particular, if $\sigma$ is a self-indexing 
Morse function, the right hand side of \eqref{eq:4.8} is just a constant.

\item If the edge, $e$, is not the only path joining $p$ to $q$, 
the identity \eqref{eq:4.7} is still true; however the right hand side of 
\eqref{eq:4.7} is in $Q(\fg_{\xi}^*)$ and has to be interpreted as the 
formal analogue of the local intersection number \eqref{eq:1.5}.

\end{enumerate}

We now return to the general case.

Let $p \stackrel{\gamma'}{\to} q$ be an ascending path from $p$ to $q$, 
let $q \stackrel{\gamma''}{\to} r$ be an ascending path from $q$ to $r$, and
let $\gamma: p \stackrel{\gamma'}{\to} q \stackrel{\gamma''}{\to} r$ be the 
ascending path from $p$ to $r$ obtained by joining $\gamma'$ and $\gamma''$. 
A direct computation shows that 
\begin{equation}
\label{eq:4.9}
\frac{E(\gamma)}{\nu_p} = \frac{E(\gamma')}{\nu_p} \cdot 
\frac{E(\gamma'')}{\nu_q} \cdot 
\frac{\alpha_{e_i}}{\rho_{e_a}(\alpha_{e_i})} \; ,
\end{equation}
where $e_i$ is the last edge of $\gamma'$ and $e_a$ is the first 
edge of $\gamma''$, both pointing upward. 

Let $\gamma : p=p_0 \to p_1 \to ... \to p_{m-1} \to p_m = q$ be an ascending 
path. We will express the contribution $E(\gamma)$ by breaking up 
the path $\gamma$ into its constituent edges. Then
\begin{equation*}
\begin{split}
\frac{E(\gamma)}{\nu_p} & = \frac{E(pp_1)}{\nu_p} \cdot ... \cdot
\frac{E(p_{m-1}q)}{\nu_{p_{m-1}}} \cdot 
\prod_{k=1}^{m-1} 
\frac{\alpha_{p_{k-1}p_k}}{\rho_{p_kp_{k+1}}(\alpha_{p_{k-1}p_k})} \\
& = \frac{-\nu_{p_1}\Theta_{pp_1}}{\alpha_{pp_1}\nu_p} \cdot
\frac{-\nu_{p_2} \Theta_{p_1p_2}}{\alpha_{p_1p_2}\nu_{p_1}} 
\cdot ... \cdot  
\frac{-\nu_{q}\Theta_{p_{m-1}q}}{\alpha_{p_{m-1}q}\nu_{p_{m-1}}} 
\prod_{k=1}^{m-1}
\frac{\alpha_{p_{k-1}p_k}}{\rho_{p_kp_{k+1}}(\alpha_{p_{k-1}p_k})} \; .
\end{split} 
\end{equation*}
Therefore the contribution $E(\gamma)$ of the path $\gamma$ is
\begin{equation}
\label{eq:4.10}
E(\gamma) = \nu_q \cdot \Bigl( \prod_{k=1}^m \Theta_{p_{k-1}p_k} \Bigr) 
\cdot \frac{(-1)^m}
{\alpha_{p_{m-1}q} \prod_{k=1}^{m-1} \rho_{p_kp_{k+1}}(\alpha_{p_{k-1}p_k})} 
\; . 
\end{equation}

In view of \eqref{eq:4.8} we can also write this in the form 
\eqref{eq:1.12}, $e_i$ being the edge of $\Gamma$ joining $p_{i-1}$ to 
$p_i$ and $\iota_e$ being the local intersection number \eqref{eq:4.7}.

If we reverse the orientation of $\Gamma$ replacing $\xi$ with $-\xi$  and the
Morse function $\phi$ by $-\phi$, we get a formula similar to \eqref{eq:1.11}
for $\tau_p^-$
\begin{equation}
\label{eq:4.11}
\tau_p^- (q) = \sum E(\gamma) \quad ,
\end{equation}
the sum being over \emph{descending} paths from $p$ to $q$.

Moreover, the $E(\gamma)$'s in \eqref{eq:4.11} are easy to compute in terms 
of the $E(\gamma)$'s in \eqref{eq:1.12}. To see this lets consider as above 
the simplest example of an ascending path in $\Gamma$, an ascending edge, 
$e$, joining $p$ to $q$. By \eqref{eq:4.3} and \eqref{eq:4.4}
\begin{align}
\theta_e E_{p,q} & = \{ e'' \in E_q^+ \; , \;  
\theta_{\bar{e}}e'' \in E_p^- \}  \NR
\intertext{and}
\theta_{\bar{e}} E_{q,p} & = \{ e' \in E_p^+ \; , \; 
\theta_e e' \in E_q^- \} - \{ e \}  \nonumber 
\end{align}
so, by \eqref{eq:4.5} and \eqref{eq:4.6}
\begin{equation}
\label{eq:4.12}
\Theta_{q,p} = \Theta_{p,q} \; . 
\end{equation}

Now let $\gamma$ be an ascending path of length $m$ from $p$ to $q$ and let 
$\bar{\gamma}$ be the same path traced in the reverse direction. Then by 
\eqref{eq:4.10} and \eqref{eq:4.12}
\begin{align}
E(\bar{\gamma})& = (-1)^m \frac{\hat{\alpha}_m}{\hat{\alpha}_1} \cdot 
\frac{\nu_p^-}{\nu_q^+} \cdot E(\gamma) \NR
\intertext{where}
\nu_p^- & = \prod_{e' \in E_p^+} \alpha_{e'}\; , \quad 
\hat{\alpha}_m = \frac{\alpha_{e_m}}{\alpha_{e_m}(\xi)} \; , \quad
\hat{\alpha}_1 = \frac{\alpha_{e_1}}{\alpha_{e_1}(\xi)}\; . \nonumber
\end{align}

\begin{figure}[h]
\begin{center}
\includegraphics{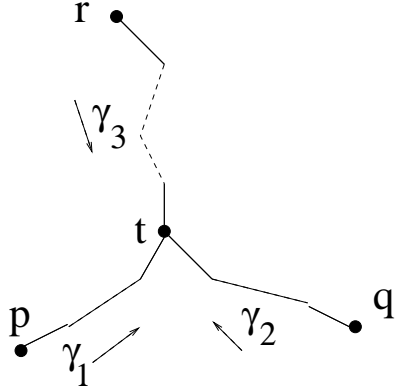}
\caption{Configuration of paths}
\label{fig:configuration}
\end{center}
\end{figure}

We are now finally in position to compute the cohomology pairing 
\eqref{eq:1.2}. By \eqref{eq:1.5}, \eqref{eq:3.12} and \eqref{eq:4.11} the 
integral
$$c_{pqr} = \int_{\Gamma} \tau_p^+ \tau_q^+ \tau_r^- $$
is equal to the sum 
$$\sum \delta_t E(\gamma_1)E(\gamma_2) E(\gamma_3)$$
summed over all triples $\gamma_1$, $\gamma_2$, $\gamma_3$ consisting of 
an ascending path, $\gamma_1$, from $p$ to $t$, an ascending path, 
$\gamma_2$, from $q$ to $t$, and a descending path, $\gamma_3$, from 
$r$ to $t$. (See Figure~\ref{fig:configuration}.) 

Thus, in particular, if there exist no such configurations, 
$c_{pqr}=0$. Now suppose that the hypothesis of Theorem~\ref{th:1.6} is 
satisfied, \emph{i.e.} $\sigma$ is a self-indexing Morse function. Then we 
claim that
\begin{equation}
\label{eq:4.13}
\int \tau_p^+ \tau_q^- = \delta_{pq} \; . 
\end{equation}
In fact if $q \not\in F_p$ the supports of $\tau_p^+$ and $\tau_q^-$ are
non-overlapping so \eqref{eq:4.13} is automatically zero; and if $q=p$,
then the support of $\tau_p^+\tau_q^-$ consists of the single point $p$
and it is easy to verify that \eqref{eq:4.13} is equal to one. Thus 
\eqref{eq:4.13} is trivially true except when $q \in F_p$ and $q \neq p$. 
In this case however, $\sigma_q > \sigma_p$ so
$$ k = \text{degree } \tau_p^+ \tau_q^- = \text{degree } \tau_p^+ +  
\text{degree } \tau_q^- = d - \sigma_q + \sigma_p < d $$
so the integral \eqref{eq:4.13} is zero just by degree considerations. 
Thus if we substitute the sum 
$$\sum c_{pq}^s \tau_s^+ $$
for $\tau_p^+ \tau_q^+$ in \eqref{eq:1.2} we obtain for $c_{pq}^s$ the 
formula \eqref{eq:1.13}.

\section{Examples}\label{sec:5}

Each of the summands in \eqref{eq:1.12} is a rational function: an element 
of the quotient field, $Q(\fg^*)$; however, the sum itself is a polynomial, 
so the singularities in the individual summands are mysteriously cancelling 
each other out. We will discuss below a few simple examples in which one can 
see how some of these cancellations are happening.

\subsection{Cancellations occuring in the individual terms} \label{ssec:5.1}

Suppose $\gamma$ is a longest ascending path from $p$ to $q$. Let 
$e_1,.. ,e_m$ be the edges of $\gamma$ ordered so that 
$t(e_{k-1}) = p_k = i(e_k)$. Then $e_k$ is the \emph{only} path joining 
$p_k$ to $p_{k+1}$; hence the intersection numbers, $\iota_{e_k}$, are 
all global intersection numbers of the form \eqref{eq:4.8} and are in 
$\SS(\fg^*)$. Hence the factor 
$$\prod \iota_{e_k} $$
in the formula \eqref{eq:1.12} is in $\SS(\fg^*)$. If, in addition, the 
Morse function $\phi$ is self-indexing, this factor is a polynomial of 
degree zero, \emph{i.e.} is just a constant.

\subsection{Nearby paths} \label{ssec:5.2}

\begin{figure}[h]
\begin{center}
\includegraphics{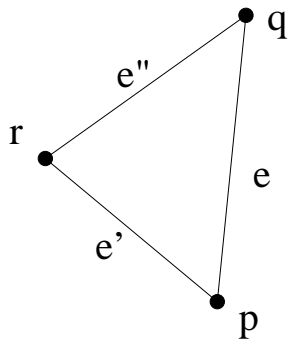}
\caption{Nearby paths}
\label{fig:nearby}
\end{center}
\end{figure}

Suppose $\Gamma$ contains a totally geodesic subgraph of the form shown in
Figure~\ref{fig:nearby}.

Let $\gamma$ be the path consisting of the single edge, $e$, joining 
$p$ to $q$ and let $\gamma_1$ be the path $p \to r \to q$. Assume $\gamma_1$
is a longest path joining $p$ to $q$ and that $\sigma_q = \sigma_p+2$. 
We claim that
$$E(\gamma) + E(\gamma_1) = \frac{\nu_q}{\alpha_e \alpha_{e''}} \; . $$

\begin{proof}
We first note that
\begin{equation}
\label{eq:5.1}
\alpha_e = \alpha_{e'} + \alpha_{e''} \; . 
\end{equation}
(This is a consequence of the compatibility conditions
$$ -\alpha_{e''} = \alpha_{e'} + c_1 \alpha_e \quad \text{and} \quad
\alpha_{e''} = \alpha_e + c_2 \alpha_{e'}$$
from which one concludes that $c_1 = c_2 = -1$.)

Let $\gamma'$ be the path joining $p$ to $r$ and $\gamma''$ the path 
joining $r$ to $q$. By \eqref{eq:4.2}
$$
E(\gamma)=\frac{\nu_q}{\alpha_e} \cdot \frac{1}{\rho_e(\alpha_{e''})} 
\; , \quad
E(\gamma') = - \frac{\nu_r}{\alpha_{e'}} \;,  \quad
E(\gamma'') = - \frac{\nu_q}{\alpha_{e''}}$$
and by \eqref{eq:4.9}
$$E(\gamma_1) = \frac{E(\gamma')}{\nu_r}\cdot 
\frac{E(\gamma'')}{\rho_{e''}(\alpha_{e'})} \alpha_{e'} \; .$$

Hence
$$E(\gamma_1) = \frac{1}{\nu_r} \Bigl( - \frac{\nu_r}{\alpha_{e'}} \Bigr) 
\Bigl( - \frac{\nu_q}{\alpha_{e''}} \Bigr) 
\frac{\alpha_{e'}}{\rho_{e''}(\alpha_{e'})} = 
\frac{\nu_q}{\alpha_{e''}\rho_{e''}(\alpha_{e'})} \; . $$
However, by \eqref{eq:5.1}, 
$\rho_{e''}(\alpha_{e'}) = \rho_{e''}(\alpha_e - \alpha_{e''}) =  
\rho_{e''}(\alpha_{e})$; hence we can rewrite this as
$$E(\gamma_1) = \frac{\nu_q}{\alpha_{e''}\rho_{e''}(\alpha_e)} \; ; $$
so $E(\gamma) + E(\gamma_1)$ is equal to the expression :
$$\nu_q \Bigl( \frac{1}{\alpha_e \rho_e(\alpha_e'')} + 
\frac{1}{\alpha_{e''} \rho_{e''}(\alpha_e)} \Bigr) \; . $$
However,
$$\rho_e(\alpha_e'') = - \frac{\alpha_{e''}(\xi)}{\alpha_e(\xi)} 
\rho_{e''}(\alpha_e ) \; , $$
so the term in parentheses can be rewritten
\begin{align}
\frac{1}{\alpha_e \rho_e(\alpha_{e''})} - &
\frac{\alpha_{e''}(\xi)}{\alpha_e(\xi)} \cdot 
\frac{1}{\alpha_{e''}\rho_e(\alpha_{e''})} = 
\frac{1}{\rho_e(\alpha_{e''})} \Bigl( \frac{1}{\alpha_e} - 
\frac{\alpha_{e''}(\xi)}{\alpha_e(\xi)} \frac{1}{\alpha_e} \Bigr) = \NR
 = & \frac{1}{\rho_e(\alpha_{e''})}  \cdot 
\frac{\alpha_{e''} - (\alpha_{e''}(\xi) / \alpha_e(\xi)) \alpha_e }
{\alpha_e \alpha_{e''}} = \frac{1}{\alpha_e \alpha_e''} \; .  \nonumber \qed 
\end{align} 
\renewcommand{\qed}{}
\end{proof}

\subsection{The flag variety $G=SL(n,\CC)/B$. } \label{ssec:5.3}

Graph theoretically, the flag variety $SL(n,\CC)/B$ is the permutahedron: 
a Cayley graph associated with the Weyl group of $SL(n,\CC)$,
the symmetric group $S_n$. Each vertex of this graph 
corresponds to a permutation $\pi \in S_n$, and two permutations 
$\pi$ and $\pi'$, are adjacent in $\Gamma$ if and only if there 
exists a transposition $\tau_{ij}$, $1 \leq i < j \leq n$ 
with $\pi'=\pi \tau_{ij}$. Moreover, if $e$ is the edge joining 
$\pi$ to $\pi\tau_{ij}$, the weight labeling $e$ is 
$$\alpha_e = \begin{cases} 
\epsilon_j - \epsilon_i ,& \mbox{ if } \pi(j)  > \pi(i) \\
\epsilon_i - \epsilon_j ,& \mbox{ if } \pi(j)  < \pi(i),
\end{cases} $$
where $\epsilon_1,..,\epsilon_n$ 
is the standard basis vectors of the lattice $\ZZ^n$.
The connection $\theta_e$ along this edge is given by
$$\theta_{\pi, \pi\tau}(\pi, \pi\tau') = (\pi \tau, \pi \tau'\tau) \; . $$

If $\xi = (\xi_1, ..., \xi_n) \in \P$, 
with $\xi_1 < ... , \xi_n$, then the function 
$$\phi : V_{\Gamma} \to \ZZ \; , \qquad \phi(\pi) = \mbox{length}(\pi)$$
is a self-indexing $\xi$-compatible Morse function on $\Gamma$.

\begin{figure}[h]
\begin{center}
\includegraphics{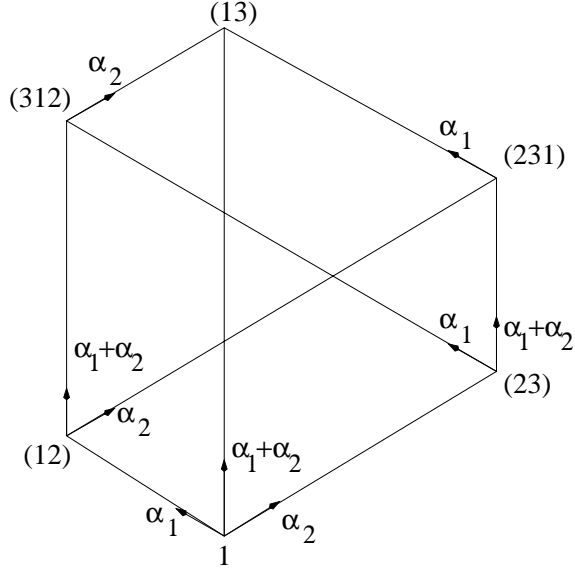}
\caption{The flag variety}
\label{fig:flag}
\end{center}
\end{figure}

The permutahedron is a bi-partite graph, with the two sets of vertices 
corresponding to even, respective odd permutations. In the special case 
$n=3$, this graph is a \emph{complete} bi-partite graph, and the 
corresponding labeling is shown in Figure~\ref{fig:flag}. 

Here $\alpha_1 = \epsilon_2-\epsilon_1$ and 
$\alpha_2 = \epsilon_3 - \epsilon_2$, and we have used the notation 
$(231)$ for the cycle $1 \to 2 \to 3 \to 1$.

The quantities $\Theta_{pq}$ given by \eqref{eq:4.6} are all equal to 1, 
with the exception of $\Theta_{1, (13)}$, which is
$$\Theta_{1,(13)} = \frac{1}{\rho_{\alpha_1+\alpha_2}(\alpha_1\alpha_2)} =
- \frac{(\alpha_1(\xi) + \alpha_2(\xi))^2}
{(\alpha_2(\xi)\alpha_1 - \alpha_1(\xi)\alpha_2)^2}$$

There are two ascending paths from (12) to (13), namely
$$\gamma_1  : (12) \to (231) \to (13) \quad \text{ and } \quad
\gamma_2  : (12) \to (312) \to (13)$$
and their contributions to $\tau_{(12)}(13)$ are 
\begin{align}
E(\gamma_1) & = - \alpha_1\alpha_2(\alpha_1+\alpha_2) 
\cdot \frac{1}{\alpha_1} 
\cdot \frac{1}{\rho_{\alpha_1}(\alpha_1+\alpha_2)} = 
\frac{\alpha_1(\xi) \alpha_2(\alpha_1 + \alpha_2)}
{\alpha_2(\xi) \alpha_1 - \alpha_1(\xi) \alpha_2} \NR
\intertext{and}
E(\gamma_2) & = -\alpha_1\alpha_2(\alpha_1+\alpha_2) 
\cdot \frac{1}{\alpha_2}
\cdot \frac{1}{\rho_{\alpha_2}(\alpha_1+\alpha_2)} = 
- \frac{\alpha_2(\xi) \alpha_1(\alpha_1 + \alpha_2)}
{\alpha_2(\xi) \alpha_1 - \alpha_1(\xi) \alpha_2} \nonumber \; ,
\end{align}
so 
$$\tau_{(12)}(13) = E(\gamma_1) + E(\gamma_2) = -\alpha_1 - \alpha_2 \; . $$

The other classes can be computed similarly and are given by
\begin{center}
\begin{tabular}{|c||c|c|c|c|c|c|} \hline
  & $\tau_1$ & $\tau_{(12)}$ & $\tau_{(23)}$ & $\tau_{(231)}$ & 
$\tau_{(312)}$ & $\tau_{(13)}$ \\ \hline \hline

 1 & 1 & 0 & 0 & 0 & 0 & 0 \\ \hline

 (12) & 1 & $-\alpha_1$ & 0 & 0 & 0 & 0 \\ \hline

 (23) & 1 & 0 & $-\alpha_2$ & 0 & 0 & 0 \\ \hline

 (231) & 1 & $-\alpha_1 - \alpha_2$ & $- \alpha_2$ &
 $\alpha_2 (\alpha_1 +\alpha_2)$ & 0 &
 0 \\ \hline

 (312) & 1 & $ -\alpha_1$ & $-\alpha_1 - \alpha_2$ &
 0 & $\alpha_1(\alpha_1 +\alpha_2)$ & 0 \\ \hline

 (13) & 1 & $-\alpha_1 -\alpha_2$ & 
 $- \alpha_1 - \alpha_2$ & $\alpha_2(\alpha_1 +\alpha_2)$ &
 $\alpha_1(\alpha_1 +\alpha_2)$ & $- \alpha_1 \alpha_2(\alpha_1 +\alpha_2)$
 \\ \hline
\end{tabular}
\end{center}

\smallskip

\subsection{The zero-dimensional Thom class} \label{ssec:5.4}

Suppose the graph $\Gamma$ is connected and hence has a unique vertex, $p_0$, 
of index zero. Then the Thom class, $\tau_{p_0}$, is the unique generator of 
$H^0(\Gamma, \alpha)$, with $\tau_{p_0}(p_0) = 1$. Thus
\begin{equation}
\label{eq:5.2}
\tau_{p_0}(p)= 1
\end{equation}
for all vertices $p$. We will show how to deduce \eqref{eq:5.2} from 
\eqref{eq:1.11}. Choose the constants, $c$ and $c'$, in \eqref{eq:3.21},
so that $p_0$ is the only vertex with $\phi (p_0) < c$ and such that 
$\phi(p)$ is the smallest critical value of $\phi$ greater that $c'$. 
By the Markov property of the map \eqref{eq:3.22}
$$1 = \sum_{w \in V_c} Q(v,w)$$
for every vertex $v \in V_{c'}$. In particular let
$$E_p^- = \{ e_i \; , \; i=1,..,r \} $$
and let $v_i \in V_{c'}$ be the vertex at which $e_i$ intersects 
$\Gamma_{c'}$. Then by \eqref{eq:3.28}
$$\tau_{p_0}(p) = \sum_{i=1}^r Q(e_i)\sum_{w}Q(v_i, w) = \sum Q(e_i)$$
and by \eqref{eq:3.26}
$$\tau_{p_0}(p) = \sum_{i=1}^r \prod_{j \neq i} 
\frac{\alpha_{e_j}}{\alpha_{e_j} - 
(\alpha_{e_j}(\xi)/\alpha_{e_i(\xi)})\alpha_{e_i}} \; . $$
Letting
$$x_i = -\frac{1}{\alpha_{e_i}(\xi)} \alpha_{e_i}$$
this becomes
$$\sum_{i=1}^r \prod_{j \neq i} \frac{-x_j}{x_i - x_j},\; , $$
which is equal to 1 by \eqref{eq:3.11}. Thus $\tau_{p_0}(p) = 1$.

\subsection{The $(n-1)$-dimensional projective space} \label{ssec:5.5}

Graph theoretically this is just the complete graph, $\Delta$, on $n$ 
vertices. Let us denote these vertices by $p_1, .., p_n$ and as in 
Section~\ref{sec:3} assign to the edge, $e$, joining $p_i$ to $p_j$, the 
weight
$$\alpha_e = x_i - x_j \; .$$
(As we pointed out in Section~\ref{sec:3} this defines an axial function on 
$\Delta$.) Let $\xi$ be an $n$-tuple of real numbers with 
$\xi_1 > \xi_2 > ... > \xi_n$ and orient the edges of $\Delta$ by decreeing 
that an edge, $e$, is ascending if $\alpha_e(\xi) > 0$. With this orientation, 
the function mapping $p_i$ to $i$ is a $\xi$-compatible Morse function. Lets 
compute the Thom class, $\tau_{p_i}$. If $i=1$, we get from the computation 
above
$$\tau_{p_1}(p) =1$$
for all vertices $p$. If $i > 1$, we can regard the vertices 
$p_i, p_{i+1}$, ... , $p_n$ as the vertices of a complete graph, $\Delta'$,
having the same axial function as above. Consider the sum 
\begin{equation}
\label{eq:5.3}
\sum E(\gamma)
\end{equation}
over all ascending paths joining $p=p_i$ to $q=p_j$, where $j > i$. 
The individual summands can be written in the form
$$ \frac{\nu_q}{\nu_q'} E'(\gamma) \; , $$
where by \eqref{eq:1.12}
$$E'(\gamma) = (-1)^m  \nu_q' \frac{\iota_{e_1}}{\alpha_m}
  \prod^m_{k=2} \frac{\iota_{e_k}}{\alpha_{k-1}-\alpha_k} $$
and $\nu_q'$ is the product
$$\prod{'}{} \alpha_{e'} $$
over all descending edge, $e' \in E_q^-$, which join $q$ to vertices in 
$\Delta'$. Then \eqref{eq:5.3} becomes
$$\frac{\nu_q}{\nu_q'} \Bigl( \sum E'(\gamma) \Bigr)  \; .$$
However, by \eqref{eq:1.11} the expression in parentheses computes the 
zeroth Thom class of the subgraph $\Delta'$, at $q$, and hence is equal 
to one. Thus
$$\tau_{p_i}(q) = \frac{\nu_q}{\nu_q'} = \prod{'}{} \alpha_{e''} \; , $$
where $\prod'$ is the product over all 
the edges $e''$ of $\Delta$, which join $q$ to the vertices 
$p_k$, $k=1,.., i-1$.

\end{document}